\documentclass[11pt]{article}
\usepackage{amssymb,amsmath,epic,eepic}
\title{
\vspace{.5ex}
\begin{center}
\LARGE {\bf Sojourn time in $\Z^+$\\ for the Bernoulli random walk on $\Z$}
\end{center}\vspace{.21ex}}

\author{Aim\'e LACHAL
\footnote{
\mbox{Postal adress: \textsc{Institut National des Sciences Appliqu\'ees de Lyon}}
\mbox{P\^ole de Math\'ematiques/Institut Camille Jordan}
\hspace*{10\textwidth} \mbox{} \mbox{B\^atiment L\'eonard de Vinci, 20 avenue Albert Einstein}
\hspace*{10\textwidth} \mbox{} \mbox{69621 Villeurbanne Cedex, \textsc{France}}
\hspace*{10\textwidth} \mbox{} E-mail: {\tt aime.lachal@insa-lyon.fr}
\hspace*{10\textwidth} \mbox{} Web page:
{\tt http://maths.insa-lyon.fr/$\mbox{}^{\sim}$lachal}
}
\\[2ex]
\begin{small}
\textsl{Universit\'e de Lyon, CNRS}
\end{small}
\\
\begin{small}
\textsl{INSA-Lyon, ICJ, UMR5208, F-69621, France}
\end{small}
}
\date{}

\newlength{\centrage}
\divide\centrage by 2
\oddsidemargin \centrage

\newtheorem{theo}{Theorem}[section]
\newtheorem{pr}[theo]{Proposition}
\newtheorem{lm}[theo]{Lemma}
\newtheorem{co}[theo]{Corollary}
\newtheorem{de}[theo]{Definition}

\newtheorem{remk}[theo]{Remark}
\newtheorem{ex}[theo]{Example}
\newenvironment{rem}{\begin{remk}\normalfont}{\ \rule{0.5em}{0.5em}\end{remk}}
\newenvironment{exam}{\begin{ex}\normalfont}{\ \rule{0.5em}{0.5em}\end{ex}}

\newcommand{\bpr}[1]{\begin{pr}#1\end{pr}}
\newcommand{\bth}[1]{\begin{theo}#1\end{theo}}
\newcommand{\blm}[1]{\begin{lm}#1\end{lm}}
\newcommand{\bco}[1]{\begin{co}#1\end{co}}

\newcommand{\brem}[1]{\begin{rem}#1\end{rem}}

\newcommand{\beq}{\begin{equation}}
\newcommand{\eeq}{\end{equation}}
\newcommand{\beqa}{\begin{align*}}
\newcommand{\eeqa}{\end{align*}}
\newcommand{\beqan}{\begin{align}}
\newcommand{\eeqan}{\end{align}}
\newcommand{\bitem}{\begin{itemize}}
\newcommand{\eitem}{\end{itemize}}

\newcommand{\dem}{\noindent {\sc Proof. }}
\newcommand{\demun}{\noindent {\sc First proof. }}
\newcommand{\demdeux}{\noindent {\sc Second proof. }}
\newcommand{\demtrois}{\noindent {\sc Third proof. }}
\newcommand{\demlem}{\noindent {\sc Proof of Lemma~\ref{lemma}. }}
\newcommand{\qed}{
\relax\ifmmode\quad\hbox{\rlap{$\sqcap$}$\sqcup$}\else
    {\unskip\nobreak\hfil\penalty50\hskip1em\null\nobreak\hfil
    \quad\hbox{\rlap{$\sqcap$}$\sqcup$}
    \parfillskip=0em\finalhyphendemerits=0\endgraf\fi}
}
\newcommand{\fin}{\ \rule{0.5em}{0.5em}}
\newcommand{\refp}[1]{(\ref{#1})}

\numberwithin{equation}{section}

\renewcommand{\a}{\alpha}

\renewcommand{\d}{\delta}

\newcommand{\e}{\varepsilon}

\renewcommand{\r}{\rho}
\newcommand{\s}{\sigma}

\newcommand{\f}{\varphi}

\newcommand{\E}{\ensuremath{\mathbb{E}}}
\newcommand{\N}{\ensuremath{\mathbb{N}}}
\renewcommand{\P}{\ensuremath{\mathbb{P}}}
\newcommand{\R}{\ensuremath{\mathbb{R}}}
\newcommand{\Z}{\ensuremath{\mathbb{Z}}}

\newcommand{\bA}{\mathbf{A}}
\newcommand{\bT}{\mathbf{T}}
\newcommand{\cE}{\mathcal{E}}
\newcommand{\cH}{\mathcal{H}}
\newcommand{\cO}{\mathcal{O}}
\newcommand{\tB}{\tilde{B}}
\newcommand{\tG}{\tilde{G}}
\newcommand{\tS}{\tilde{S}}

\renewcommand{\ge}{\geqslant}
\renewcommand{\le}{\leqslant}
\newcommand{\dis}{\displaystyle}

\newcommand{\ind}{1\hspace{-.27em}\mbox{\rm l}}
\newcommand{\lqn}[1]{\noalign{\noindent $\displaystyle{#1}$}}
\newcommand{\sN}{{\scriptscriptstyle N}}
\newcommand{\eN}{\e_{_N}}
\newcommand{\pN}{p_{_N}}
\newcommand{\qN}{q_{_N}}

\begin{document}
\maketitle

%%%%%%%%%%%%%%%%%%%%%%%%%%%%%%%%%%%%%%%%%%%%%%%%%%%%%%%%%%%%%%%%%%%%%%%%%%%%%%%
%%%%%%%%%%%%%%%%%%%%%%%%% RESUME %%%%%%%%%%%%%%%%%%%%%%%%%%%%%%%%%%%%%%%%%%%%%%
%%%%%%%%%%%%%%%%%%%%%%%%%%%%%%%%%%%%%%%%%%%%%%%%%%%%%%%%%%%%%%%%%%%%%%%%%%%%%%%
\begin{abstract}
Let $(S_k)_{k\ge 1}$ be the classical Bernoulli random walk on the integer line
with jump parameters $p\in(0,1)$ and $q=1-p$. The probability distribution of
the sojourn time of the walk in the set of non-negative integers up to a fixed
time is well-known, but its expression is not simple. By modifying slightly this
sojourn time--through a particular counting process of the zeros of the
walk as done by Chung \& Feller [``On fluctuations in coin-tossings'',
Proc. Nat. Acad. Sci. U.S.A. 35 (1949), 605--608]--, simpler
representations may be obtained for its probability distribution.
In the aforementioned article, only the symmetric case ($p=q=1/2$) is considered.
This is the discrete counterpart to the famous Paul L\'evy's arcsine law for
Brownian motion.

In the present paper, we write out a representation for this probability
distribution in the general case together with others related to the random
walk subject to a possible conditioning.
The main tool is the use of generating functions.
\end{abstract}

\begin{footnotesize}\sc
\noindent AMS 2000 subject classifications:
\begin{tabular}[t]{l}
primary 60G50; 60J22; \\ secondary 60J10; 60E10.
\end{tabular}  
\end{footnotesize}

%%%%%%%%%%%%%%%%%%%%%%%%%%%%%%%%%%%%%%%%%%%%%%%%%%%%%%%%%%%%%%%%%%%%%%%%%%%%%%%
%%%%%%%%%%%%%%%%%%%%%%%%%    Introduction      %%%%%%%%%%%%%%%%%%%%%%%%%%%%%%%%
%%%%%%%%%%%%%%%%%%%%%%%%%%%%%%%%%%%%%%%%%%%%%%%%%%%%%%%%%%%%%%%%%%%%%%%%%%%%%%%
\newpage
\tableofcontents

\section{Introduction}

Let $(X_k)_{k\ge 1}$ be a sequence of Bernoulli random variables with parameters
$p=\mbox{$\P\{X_k=1\}$}\in(0,1)$ and $q=1-p=\P\{X_k=-1\}$, and $(S_k)_{k\ge 0}$
be the random walk defined on the set of integers $\Z=\{\dots,-1,0,1,\dots\}$
as $S_k=S_0+\sum_{j=1}^k X_j$, $k\ge 1$ with initial location $S_0$.
For brevity, we write $\P_i=\P\{\dots|S_0=i\}$ and $\P_0=\P$.

The probability distribution of the sojourn time of the walk $(S_k)_{k\ge 0}$
in $\Z^+=\Z\cap [0,+\infty)$ up to a fixed step $n\ge 1$,
$N_n=\sum_{j=0}^n \ind_{\Z^+}(S_j)=\#\{j\in\{0,\dots,n\}:S_j\ge 0\}$,
is well-known.
A representation for this probability distribution can be derived with
the aid of Sparre Andersen's theorem (see~\cite{sparre1,sparre2} and,
e.g., \cite[Chap. IV, \S20]{spitzer}).
This latter can be stated according as the remarkable relationship, setting $N_0=0$,
$$
\P\{N_n=k\}=\P\{N_k=k\}\P\{N_{n-k}=0\} \mbox{ for } 0\le k \le n
$$
where the probabilities $\P\{N_k=0\}$ and $\P\{N_k=k\}$, $k\in\N$, are implicitly
known through their generating functions:
\begin{align*}
\sum_{k=0}^{\infty}\P\{N_k=0\}z^k
&=
\exp\bigg[\sum_{k=1}^{\infty}\P\{S_k<0\}\,\frac{z^k}{k}\bigg],
\\
\sum_{k=0}^{\infty}\P\{N_k=k\}z^k
&=
\exp\bigg[\sum_{k=1}^{\infty}\P\{S_k\ge 0\}\,\frac{z^k}{k}\bigg].
\end{align*}
Nevertheless, the result is not so simple.
Rescaling the random walk and passing to the limit, we get the most famous
Paul Levy's arcsine law for Brownian motion.

By modifying slightly the counting process of the positive terms of the random
walk as done by Chung \& Feller (see~\cite{chung} and, e.g., \cite[Chap. III, \S4]{feller}
and~\cite[Chap. 8, \S11]{renyi}), an alternative sojourn time of the walk
$(S_k)_{k\in\N}$ in $\Z^+$ up to $n$ can be defined as $T_n=\sum_{j=1}^n \d_j$ with
$$
\d_j=\left\{\begin{array}{ll}
1 & \mbox{if $(S_j>0)$ or $(S_j=0$ and $S_{j-1}>0)$,}
\\
0 & \mbox{if $(S_j<0)$ or $(S_j=0$ and $S_{j-1}<0)$.}
\end{array}\right.
$$
We put $T_0=0$. We obviously have $0\le T_n\le n$.
In $T_n$, $n\ge 1$, one counts each step $j$ such that $S_j>0$ and only those
steps such that $S_j=0$ which correspond to a downstep: $S_{j-1}=1$.
%In this situation, a simple formula for the probability distribution of $T_n$ can be obtained.
This convention is described in~\cite{chung} (and, e.g., in~\cite{feller}
and~\cite{renyi}) in the symmetric case $p=q=1/2$ when $n$ is an even integer
and, as written in~\cite{chung}--\textit{``The elegance of the results to be
announced depends on this convention''} (\textit{sic})--, it produces a remarkable result.
Indeed, in this case, the sojourn time is even and its probability distribution takes the
simple following form: for even integers $k$ such that $0\le k\le n$,
as in Sparre Andersen's theorem,
$$
\P\{T_n=k\}=\P\{T_k=k\}\P\{T_{n-k}=0\}=\frac{1}{2^n}\binom{k}{k/2}\binom{(n-k)}{(n-k)/2}.
$$

In this paper, we derive explicit expressions for the probability distribution
of $T_n$ in any case, that is for any $p\in(0,1)$ and any integer (even or odd)
$n\ge 1$. The main results are displayed in Theorems~\ref{theorem0bis} and~\ref{theorem1bis}.
We also compute the distribution of $T_n$ under various constraints at
the last step: $S_n=0$, $S_n>0$ or $S_n<0$. The constraint $S_n=0$ (with $S_n=0$)
corresponds to the bridge of the random walk.
The related results are respectively included in Theorems~\ref{theorem2bis} and~\ref{theorem3bis}.
We examine in details several examples corresponding to the cases $n\in\{0,1,2,\dots,8\}$.
The main tool for this study is the use of generating functions together
with clever algebra. The intermediate results are contained in
Theorems~\ref{theorem0}, \ref{theorem1}, \ref{theorem2} and~\ref{theorem3}.
For several results, we find it interesting to produce two or three proofs.
Certain are direct while others rely on recursive properties.

Finally, be rescaling suitably the random walk, we retrieve the distribution
of the sojourn time in $(0,+\infty)$ for the Brownian motion with a possible
drift. This includes of course the famous Paul L\'evy's arcsine law for
undrifted Brownian motion.

Although this problem is old and classical, we are surprised not to have
found any related reference in the literature.

%%%%%%%%%%%%%%%%%%%%%%%%%%%%%%%%%%%%%%%%%%%%%%%%%%%%%%%%%%%%%%%%%%%%%%%%%%%%%%%
%%%%%%%%%%%%%%%%%%%%%%%%%        Notation      %%%%%%%%%%%%%%%%%%%%%%%%%%%%%%%%
%%%%%%%%%%%%%%%%%%%%%%%%%%%%%%%%%%%%%%%%%%%%%%%%%%%%%%%%%%%%%%%%%%%%%%%%%%%%%%%
\section{Settings and mathematical background}

%%%%%%%%%%%%%%%%%%%%%%%%%%%%%%%%%%%%%%%%%%%%%%%%%%%%%%%%%%%%%%%%%%%%%%%%%%%%%%%
\subsection{Some preliminary identities}

Let $\N=\Z^+=\{0,1,2,\dots\}$ be the usual set of non-negative integers,
$\N^*=\Z^{+*}=\N\setminus\{0\}=\{1,2,3,\dots\}$ that of positive integers,
$\Z^{-*}=\{\dots,-3,-2,-1\}$ that of negative integers and $\Z^*=\Z\setminus\{0\}$.
Let $\cE=\{0,2,4,\dots\}$ denote the set of even non-negative integers,
$\cE^*=\cE\setminus\{0\}=\{2,4,6,\dots\}$ the set of even positive integers and
$\cO=\{1,3,5,\dots\}$ the set of odd positive integers. Set, for suitable real $z$,
$$
A(z)=\sqrt{1-4z^2} \quad\mbox{and}\quad \bA(z)=\sqrt{1-4pqz^2}.
$$
Set also for $i\in\cE$
$$
\dis a_i=\frac{1}{i+2}\binom{i}{i/2}=\frac{1}{4(i+1)}\binom{i+2}{(i+2)/2}
\quad\mbox{and}\quad b_i=(i+2)a_i=\binom{i}{i/2}.
$$
We shall make use of the following elementary identities.
%
%%%%%%%%%%%%%%
\bpr{
For any $z$ such that $|z|<1/2$,
\begin{align}
A(z)
&=
\sqrt{1-4z^2}=-\sum_{i\in\cE} \frac{1}{i-1}\binom{i}{i/2}z^i
=1-4\sum_{i\in\cE^*} a_{i-2}z^i,
\nonumber\\[-1ex]
&
\label{elem}\\[-2ex]
\frac{1}{A(z)}
&=
\frac{1}{\sqrt{1-4z^2}}=\sum_{i\in\cE} \binom{i}{i/2}z^i=\sum_{i\in\cE} b_i z^i.
\nonumber
\end{align}
}
%%%%%%%%%%%%%%
%
We have the following convolution relationships.
%
%%%%%%%%%%%%%%
\bpr{
For any even integer $i\ge 0$,
\beq
\sum_{j\in\cE:\atop j\le i}a_ja_{i-j}=\frac{1}{2}\,a_{i+2}
\quad\mbox{and}\quad \sum_{j\in\cE:\atop j\le i}a_jb_{i-j}=\frac{1}{4}\,b_{i+2}.
\label{convol}
\eeq
}
%%%%%%%%%%%%%%
%
\dem
The generating function of the left-hand side of the first equality
in~\refp{convol} can be evaluated as
\begin{align*}
\sum_{i\in\cE}\bigg(\sum_{j\in\cE:\atop j\le i}a_ja_{i-j}\bigg)z^i
&=
\bigg(\sum_{i\in\cE} a_iz^i\bigg)^2=\bigg(\frac{1-A(z)}{4z^2}\bigg)^2
\\
&=
\frac{1}{8z^4} (1-2z^2-A(z))=\frac{1}{2} \sum_{i\in\cE} a_{i+2}z^i.
\end{align*}
By identification of the coefficients of the foregoing generating functions, we
immediately obtained the first equality in~\refp{convol}.
Analogously, for the second equality in~\refp{convol},
\begin{align*}
\sum_{i\in\cE}\bigg(\sum_{j\in\cE:\atop j\le i}a_jb_{i-j}\bigg)z^i
&=
\bigg(\sum_{i\in\cE} a_iz^i\bigg)\bigg(\sum_{i\in\cE} b_iz^i\bigg)
=\frac{1-A(z)}{4z^2A(z)}
\\
&=
\frac{1}{4z^2A(z)}-\frac{1}{4z^2}=\frac{1}{4} \sum_{i\in\cE} b_{i+2}z^i
\end{align*}
and the second equality in~\refp{convol} holds.
\fin

\noindent We shall also use the identity below.
%
%%%%%%%%%%%%%%
\bpr{
For any $x$ and $y$ such that $|x|<1/2$ and $|y|<1/2$,
\beq
\frac{1}{A(x)+A(y)}=\sum_{i,j\in\cE}a_{i+j} \,x^iy^j.
\label{1surA+A}
\eeq
}
%%%%%%%%%%%%%%
%
\dem
Let us write
$\dis \frac{1}{A(x)+A(y)}=\frac{A(x)-A(y)}{4(y^2-x^2)}.$
On one hand,
$$
A(x)-A(y)=4\sum_{k\in\cE} a_k (y^{k+2}-x^{k+2}).
$$
On the other hand,
$$
\frac{y^{k+2}-x^{k+2}}{y^2-x^2}=\sum_{i,j\in\cE:\atop i+j=k} x^iy^j.
$$
Therefore,
$$
\frac{1}{A(x)+A(y)}=\sum_{k\in\cE} a_k \sum_{i,j\in\cE:\atop i+j=k} x^iy^j=
\sum_{i,j\in\cE} a_{i+j}\,x^iy^j.
$$
\fin

%%%%%%%%%%%%%%%%%%%%%%%%%%%%%%%%%%%%%%%%%%%%%%%%%%%%%%%%%%%%%%%%%%%%%%%%%%%%%%%
\subsection{Some well-known identities on random walks}

Here, we recall several well-known formulas in the theory of random walks.
We refer, e.g., to~\cite{spitzer}.
We have, for $j\in\Z$ and $k\in\N$ such that $|j|\le k$ and $k-j\in\cE$,
\beq
\P\{S_k=j\}=\binom{k}{(j+k)/2} p^{(j+k)/2}q^{(k-j)/2}.
\label{law-Sk}
\eeq
Using the representation $\binom{k}{(k+1)/2}=b_{k+1}/2$, we have in particular
\begin{align}
\P\{S_k=0\}
&=
\left\{\begin{array}{ll}
b_k (pq)^{k/2} &\mbox{if $k$ is even,}
\\[0ex]
0 &\mbox{if $k$ is odd,}
\end{array}\right.
\nonumber
\\[1ex]
\P\{S_k=1\}=\P_{-1}\{S_k=0\}
&=
\left\{\begin{array}{ll}
\frac12\,b_{k+1}\, p^{(k+1)/2}q^{(k-1)/2} &\mbox{if $k$ is odd,}
\\[1ex]
0 &\mbox{if $k$ is even,}
\end{array}\right.
\label{law-Sk-part}
\\[1ex]
\P\{S_k=-1\}=\P_1\{S_k=0\}
&=
\left\{\begin{array}{ll}
\frac12\,b_{k+1}\, p^{(k-1)/2}q^{(k+1)/2} &\mbox{if $k$ is odd,}
\\[1ex]
0 &\mbox{if $k$ is even.}
\end{array}\right.
\nonumber
\end{align}

We define several generating functions.
For $j\in\Z$, let $H_j$ be the generating function of the $\P\{S_k=j\},k\in\N$,
and, for $i\in\Z$, $H_i^F$ be the generating function of the $\P_i\{S_k\in F\},$ $k\in\N$:
\beq
H_j(z)=\sum_{k=0}^{\infty}\P\{S_k=j\}z^k \quad\mbox{and}\quad
H_i^F(z)=\sum_{k=0}^{\infty}\P_i\{S_k\in F\}z^k=\sum_{j\in F-i} H_j(z)
\label{def-fct-gene}
\eeq
where $F-i$ is the set of the numbers of the form $j-i,i\in F$. We explicitly have
\beq
H_j(z)=\left\{\begin{array}{ll}
\dis\frac{1}{\bA(z)}\left(\frac{1-\bA(z)}{2qz}\right)^j
=\frac{1}{\bA(z)}\left(\frac{2pz}{1+\bA(z)}\right)^j & \mbox{if $j\ge 0$,}
\\[2ex]
\dis\frac{1}{\bA(z)}\left(\frac{1-\bA(z)}{2pz}\right)^{|j|}
=\frac{1}{\bA(z)}\left(\frac{2qz}{1+\bA(z)}\right)^{|j|} & \mbox{if $j\le 0$.}
\end{array}\right.
\label{geneH}
\eeq

We need to introduce the first hitting time of a level $a\in\Z$ for the random
walk: $\tau_a=\min\{n\ge 1:S_n=a\}$. The probability distribution of $\tau_a$
for $a\in\Z^*$ can be expressed by means of the probabilities $\P\{S_k=a\}$,
$k\in\N$, according as
\beq
\P\{\tau_a=k\}=\frac{|a|}{k}\,\P\{S_k=a\}
=\frac{|a|}{k}\binom{k}{(k+a)/2}p^{(k+a)/2}q^{(k-a)/2} \mbox{ for $k\ge |a|$}.
\label{law-tau-a}
\eeq
In some particular cases, we have for $k\in\cE^*$:
$$
\pi_k=\P\{\tau_0=k\}=\frac{1}{k-1}\binom{k}{k/2}(pq)^{k/2}
$$
(and we set $\pi_0=0$ and $\pi_k=0$ for $k\in \cO$) and for $k\in \cO$:
\begin{align*}
\P\{\tau_1=k\}&=\frac{1}{k}\binom{k}{(k+1)/2}p^{(k+1)/2}q^{(k-1)/2},
\\
\P\{\tau_{-1}=k\}&=\frac{1}{k}\binom{k}{(k+1)/2}p^{(k-1)/2}q^{(k+1)/2}.
\end{align*}
We sum up these formulas as follows:
\bitem
\item
for $k\in\cE^*$,
\beq
\pi_k=4a_{k-2}(pq)^{k/2};
\label{law-tau0}
\eeq
\item
for $k\in \cO$,
\begin{align}
\P\{\tau_1=k\}&=\P_{-1}\{\tau_0=k\}=\frac{1}{2q}\,\pi_{k+1}=2p\,a_{k-1}(pq)^{(k-1)/2},
\nonumber\\[-1ex]
\label{law-tau1}
\\[-2ex]
\P\{\tau_{-1}=k\}&=\P_1\{\tau_0=k\}=\frac{1}{2p}\,\pi_{k+1}=2q\,a_{k-1}(pq)^{(k-1)/2}.
\nonumber
\end{align}
\eitem

%
%%%%%%%%%%%%%%
\brem{
The convolution identities~\refp{convol} can be interpreted as
Darling-Siegert-type equations (see, e.g., \cite{darling}) which are due to the
Markov property of the random walk. More precisely, the second identity
of~\refp{convol} is the analytic form of the probabilistic equality
$$
\P\{S_n=0\}=\sum_{j\in\cE:\atop j\le n} \P\{\tau_0=j\}\P\{S_{n-j}=0\}
$$
which is obtained by remarking that a trajectory starting at zero
and terminating at zero at time $n$ necessarily passes through zero at
a time equal to or less than $n$ and then $\tau_0\le n$.
The first identity of~\refp{convol} is the analytic form
of the probabilistic equality
$$
\P\{\tau_2=n\}=\sum_{j\in\cO:\atop j\le n-1} \P\{\tau_1=j\}\P\{\tau_1=n-j\}
$$
which is obtained by observing that a trajectory starting at zero and passing
at level two at time $n$ for the first time necessarily crosses level
one at a time less than $n$: $\tau_1<n$.
}
%%%%%%%%%%%%%%
%

The corresponding generating functions $K_a$ defined as
$$
K_a(z)=\E(z^{\tau_a})=\sum_{k=1}^{\infty} \P\{\tau_a=k\}z^k
$$
are explicitly expressed by
\beq
K_a(z)=\left\{\begin{array}{ll}
\dis\left(\frac{1-\bA(z)}{2qz}\right)^a & \mbox{if $a\ge 1$,}
\\[2ex]
\dis\left(\frac{1-\bA(z)}{2pz}\right)^{|a|} & \mbox{if $a\le -1$,}
\\[2ex]
1-\bA(z) & \mbox{if $a=0$.}
\end{array}\right.
\label{geneK}
\eeq
We state a last elementary result. Noticing, by the Markov property, that
$\P\{\mbox{$\tau_0=j$},\mbox{$S_1>0$}\}$ $=p\,\P_1\{\tau_0=j-1\}$ and
$\P\{\tau_0=j,S_1<0\}=q\,\P_1\{\tau_0=j-1\},$
we get, by~\refp{law-tau1},
\beq
\P\{\tau_0=j,S_1>0\}=\P\{\tau_0=j,S_1<0\}=\frac12 \,\pi_j.
\label{law-tau0-S1}
\eeq

%%%%%%%%%%%%%%%%%%%%%%%%%%%%%%%%%%%%%%%%%%%%%%%%%%%%%%%%%%%%%%%%%%%%%%%%%%%%%%%
\subsection{Purpose of the article}

The aim of this paper is to compute the probabilities
$$
r_{k,n}^F=\P\{T_n=k,S_n\in F\}, 0\le k\le n,
$$
in the four cases $F=\Z$, $F=\{0\}$, $F=\Z^{+*}$ and $F=\Z^{-*}$.
The corresponding results are displayed in Theorems~\ref{theorem1},
\ref{theorem2} and~\ref{theorem3}.
For this, we write out a recurrence relation (see Proposition~\ref{prop-rec})
from which we derive the result.
We provide two approaches: one hinges on a method by induction, the
other--which is more constructive--consists of first calculating the
generating function $G_F$ of the $r_{k,n}^F$'s:
$$
G^F(x,y)=\sum_{k,n\in\N:\atop k\le n}r_{k,n}^F \,x^ky^{n-k}
=\sum_{n\in\N} \E\left(x^{T_n} y^{n-T_n} \ind_{\{S_n\in F\}}\right)
$$
and next of inverting this function. Notice that $r_{0,0}^F=\ind_F(0)$ and for
$n\ge 1$, the events $\{T_n=0\}$ and $\{T_n=n\}$ respectively coincide with
$\{S_1\le 0,\dots,S_n\le 0\}=\{\tau_1>n\}$ and
$\{S_1\ge 0,\dots,S_n\ge 0\}=\{\tau_{-1}>n\}$. Therefore, for $n\ge 1$,
\beq
r_{0,n}^F=\P\{\tau_1>n,S_n\in F\}\quad\mbox{and}\quad
r_{n,n}^F=\P\{\tau_{-1}>n,S_n\in F\}.
\label{r0nF}
\eeq
The particular probabilities $r_{0,n}^F$ and $r_{n,n}^F$, $n\in\N$, are
generated by the partial functions of $G^F$, namely
\beq
G^F(x,0)=\sum_{n\in\N} r_{n,n}^F \,x^k\quad\mbox{and}\quad
G^F(0,y)=\sum_{n\in\N} r_{0,n}^F \,y^n.
\label{partial}
\eeq

%%%%%%%%%%%%%%%%%%%%%%%%%%%%%%%%%%%%%%%%%%%%%%%%%%%%%%%%%%%%%%%%%%%%%%%%%%%%%%%
%%%%%%%%%%%%%%%%%%%%%%%%%     Recurrence       %%%%%%%%%%%%%%%%%%%%%%%%%%%%%%%%
%%%%%%%%%%%%%%%%%%%%%%%%%%%%%%%%%%%%%%%%%%%%%%%%%%%%%%%%%%%%%%%%%%%%%%%%%%%%%%%
\section{A recurrence relationship}

As in~\cite{chung}, we state the recursive relationship below.
%
%%%%%%%%%%%%%%
\bpr{\label{prop-rec}
The following relationship holds for $1\le k\le n-1$:
\beq
r_{k,n}^F=\frac12 \bigg[\sum_{j\in\cE^*:\atop j\le k}\pi_j\,r_{k-j,n-j}^F
+\sum_{j\in\cE^*:\atop j\le n-k}\pi_j\,r_{k,n-j}^F\bigg].
\label{rec}
\eeq
For $k=0$ and $k=n$, the corresponding $r_{k,n}^F$ admit the following representations:
\begin{align}
r_{0,n}^F
&=
\P\{S_n\in F\}-\frac{1}{2q}
\sum_{j\in\cO:\atop j\le n} \pi_{j+1}\P_1\{S_{n-j}\in F\},
\nonumber\\[-2ex]
&
\label{part-cases}
\\[-1ex]
r_{n,n}^F
&=
\P\{S_n\in F\}-\frac{1}{2p}
\sum_{j\in\cO:\atop j\le n} \pi_{j+1}\P_{-1}\{S_{n-j}\in F\},
\nonumber
\end{align}
where $\pi_j=4a_{j-2}(pq)^{j/2}$.
}
%%%%%%%%%%%%%%
%
\dem
Let us consider the event $\{T_n=k\}$. If $1\le k\le n-1$, the random
walk passes through level 0 before the step $n$, that is $1\le \tau_0\le n-1$.
Then $T_n=\sum_{j=1}^{\tau_0} \d_j+T'_{k,n}$ where
$T'_{k,n}=\sum_{j=\tau_0+1}^n \d_j$.
Now, if the first step is positive ($S_1>0$), then $T_n=\tau_0+T'_{k,n}$
and if the first step is negative ($S_1<0$), then $T_n=T'_{k,n}$.
This discussion yields, with the aid of the Markov property,
\begin{align*}
r_{k,n}^F
&=
\P\{T_n=k,S_n\in F,1\le\tau_0 \le n-1\}
=\sum_{j=1}^{n-1} \P\{T_n=k,S_n\in F,\tau_0=j\}
\\
&=
\sum_{j\in\cE^*:\atop j\le k} \P\{\tau_0=j,S_1>0\} \P\{T_{n-j}=k-j,S_{n-j}\in F\}
\\
&\hphantom{=\,}
+\!\sum_{j\in\cE^*:\atop j\le n-k} \P\{\tau_0=j,S_1<0\} \P\{T_{n-j}=k,S_{n-j}\in F\}.
\end{align*}
We conclude that~\refp{rec} holds thanks to~\refp{law-tau0-S1}.

For the cases $k=0$ and $k=n$, by~\refp{r0nF}, the probabilities of
interest can be expressed as follows:
\begin{align*}
r_{0,n}^F
&=
\P\{S_n\in F\}-\P\{\tau_1\le n,S_n\in F\}
=\P\{S_n\in F\}-\sum_{j=1}^n \P\{\tau_1=j\} \P_1\{S_{n-j}\in F\},
\\
r_{n,n}^F
&=
\P\{S_n\in F\}-\P\{\tau_{-1}\le n,S_n\in F\}
=\P\{S_n\in F\}-\sum_{j=1}^n \P\{\tau_{-1}=j\} \P_{-1}\{S_{n-j}\in F\}.
\end{align*}
In view of~\refp{law-tau1}, this proves~\refp{part-cases}.
\fin
%%%%%%%%%%%%%%

The relationships~\refp{rec} and~\refp{part-cases} allow to exhibit an explicit
expression for $G^F$.
%%%%%%%%%%%%%%
\bth{\label{th-generalG}
The generating function $G^F$ can be written as
\beq
G^F(x,y)=\frac{[1+\bA(x)]G^F(x,0)+[1+\bA(y)]G^F(0,y)-2\ind_F(0)}{\bA(x)+\bA(y)}
\label{generalG}
\eeq
where $\bA(z)=\sqrt{1-4pqz^2}$.
The quantities $G^F(x,0)$ and $G^F(0,y)$ can be expressed as
\beq
G^F(x,0)=H_0^F(x)-\frac{1-\bA(x)}{2px}\,H_{-1}^F(x)\quad\mbox{and}\quad
G^F(0,y)=H_0^F(y)-\frac{1-\bA(y)}{2qy}\,H_1^F(x)
\label{GF(x,0)}
\eeq
where $H_0^F$, $H_1^F$ and $H_{-1}^F$ are defined by~\refp{def-fct-gene}.
}
%%%%%%%%%%%%%%
\dem
We decompose the sum defining $G^F$ into four parts:
\begin{align}
G^F(x,y)
&=
\sum_{n\in\N} r_{0,n}^F y^n+\sum_{n\in\N} r_{n,n}^Fx^n -r_{0,0}^F
+\sum_{k,n\in\N:\atop 1\le k\le n-1} r_{k,n}^F x^ky^{n-k}
\nonumber\\
&=
G^F(x,0)+G^F(0,y)-\ind_F(0)+\sum_{k,n\in\N:\atop 1\le k\le n-1}
r_{k,n}^F x^ky^{n-k}.
\label{cut-GF}
\end{align}
In view of~\refp{rec}, the sum in~\refp{cut-GF} can be written as
\beq
\sum_{k,n\in\N:\atop 1\le k\le n-1} r_{k,n}^F x^ky^{n-k}
=\frac12 \bigg[\sum_{j,k,n\in\cE^*:\atop 1\le j\le k\le n-1}
\pi_j\,r_{k-j,n-j}^F x^ky^{n-k}
+\sum_{j,k,n\in\cE^*:\atop 1\le j\le n-k}\pi_j\,r_{k,n-j}^F x^ky^{n-k}\bigg].
\label{sum1}
\eeq
We have
\begin{align}
\sum_{j,k,n\in\cE^*:\atop 1\le j\le k\le n-1}\pi_j\,r_{k-j,n-j}^F x^ky^{n-k}
&=
\sum_{j\in\cE^*}\pi_j\,x^j\sum_{k,n\in\cE^*:\atop j\le k\le n-1}
r_{k-j,n-j}^F x^{k-j}y^{(n-j)-(k-j)}
\nonumber\\
&=
\sum_{j\in\cE^*}\pi_j\,x^j \sum_{k,n\in\cE^*:\atop k\le n-1}
r_{k,n}^F x^{k}y^{n-k}
\nonumber\\
&=
K_0(x)[G^F(x,y)-G^F(x,0)]
\label{sum2}
\end{align}
where $K_0$ is defined in~\refp{geneK}. Similarly,
\begin{align}
\sum_{j,k,n\in\cE^*:\atop 1\le j\le n-k}\pi_j\, r_{k,n-j}^F x^ky^{n-k}
&=
\sum_{j\in\cE^*}\pi_j\,y^j\sum_{k,n\in\cE^*:\atop k\le n-j}
r_{k,n-j}^F x^{k}y^{(n-j)-k}
\nonumber\\
&=
\sum_{j\in\cE^*}\pi_j\, y^j \sum_{k,n\in\cE^*:\atop k\le n} r_{k,n}^F x^{k}y^{n-k}
\nonumber\\
&=
K_0(y)[G^F(x,y)-G^F(0,y)].
\label{sum3}
\end{align}
By plugging~\refp{sum2} and~\refp{sum3} into~\refp{sum1} and next
into~\refp{cut-GF}, we obtain the equation for $G^F(x,y)$
$$
G^F(x,y)=[1-\frac12\,K_0(x)]G^F(x,0)+[1-\frac12\,K_0(y)]G^F(0,y)-\ind_F(0)
+\frac12\,[K_0(x)+K_0(y)]G^F(x,y)
$$
from which we immediately extract the result~\refp{generalG}.

Concerning the terms $G^F(x,0)$ and $G^F(0,y)$, we have
by~\refp{partial} and~\refp{part-cases}
$$
G^F(x,0)=\sum_{n\in\N} r_{n,n}^Fx^n
=\sum_{n\in\N}\P\{S_n\in F\}x^n-\frac{1}{2p}\sum_{n\in\N^*}
\sum_{j\in\cO:\atop j\le n} \pi_{j+1} \P_{-1}\{S_{n-j}\in F\}x^n.
$$
The first sum in the above equality is $H_0^F(x)$ while the second can
be computed as follows: by~\refp{law-tau1},
\begin{align*}
\sum_{n\in\N^*}\sum_{j\in\cO:\atop j\le n} \pi_{j+1} \P_{-1}\{S_{n-j}\in F\}x^n
&=
\sum_{j\in\cO} \pi_{j+1} x^j \sum_{n\in\N^*:\atop n\ge j}
\P_{-1}\{S_{n-j}\in F\}x^{n-j}
\\
&=
2p \sum_{j\in\cO}\P\{\tau_{-1}=j\} x^j \sum_{n\in\N} \P_{-1}\{S_n\in F\}x^n
\\
&=
2p K_{-1}(x)H_{-1}^F(x)
\end{align*}
and the expression of $G^F(x,0)$ in~\refp{GF(x,0)} ensues. The derivation of
that of $G^F(0,y)$ is quite similar.
The proof of Theorem~\ref{th-generalG} is finished.
\fin

%%%%%%%%%%%%%%%%%%%%%%%%%%%%%%%%%%%%%%%%%%%%%%%%%%%%%%%%%%%%%%%%%%%%%%%%%%%%%%%
%%%%%%%%%%%%%%%%%%%%%  Pinned random walk      %%%%%%%%%%%%%%%%%%%%%%%%%%%%%%%%
%%%%%%%%%%%%%%%%%%%%%%%%%%%%%%%%%%%%%%%%%%%%%%%%%%%%%%%%%%%%%%%%%%%%%%%%%%%%%%%
\section{The case $F=\{j\}$ for $j\in\Z^*$}

We suppose that $F=\{j\}$ for a fixed $j\in\Z^*$. So, we are dealing
with a random walk with a prescribed location after the $n$th step.
The case where $j=0$ will be
considered in Section~\ref{section-bridge}. We set for simplicity
$r_{k,n}^{\{j\}}=r_{k,n}^j$ and $G^{\{j\}}(x,y)=G^j(x,y)$,
$H_i^{\{j\}}(z)=H_i^j(z)$.

%%%%%%%%%%%%%%%%%%%%%%%%%%%%%%%%%%%%%%%%%%%%%%%%%%%%%%%%%%%%%%%%%%%%%%%%%%%%%%%
%%%%%%%%%%%%%%%%%%%%%%%%%  Generating function   %%%%%%%%%%%%%%%%%%%%%%%%%%%%%%
%%%%%%%%%%%%%%%%%%%%%%%%%%%%%%%%%%%%%%%%%%%%%%%%%%%%%%%%%%%%%%%%%%%%%%%%%%%%%%%
\subsection{Generating function}

In order to write the generating function $G^j(x,y)$, in view of~\refp{generalG},
we need to know that $H_i^j(z)=H_{j-i}(z)$ and to evaluate the functions
$G^j(x,0)$ and $G^j(0,y)$. We have
$$
G^j(x,0)=H_0^j(x)-\frac{1-\bA(x)}{2px}\,H_{-1}^j(x).
$$
But
$$
H_{-1}^j(x)=\left\{\begin{array}{ll}
\dis \frac{2px}{1+\bA(x)}\,H_0^j(x) & \mbox{if $j\in\Z^{+*}$,}
\\[3ex]
\dis \frac{2px}{1-\bA(x)}\,H_0^j(x) & \mbox{if $j\in\Z^{-*}$.}
\end{array}\right.
$$
Then, for $j\in\Z^{+*}$,
$$
G^j(x,0)=H_0^j(x)\bigg[1-\frac{1-\bA(x)}{2px}\,\frac{2px}{1+\bA(x)}\bigg]
=\frac{2K_j(x)}{1+\bA(x)}
$$
and, for $j\in\Z^{+*}$,
$$
G^j(x,0)=H_0^j(x)\bigg[1-\frac{1-\bA(x)}{2px}\,\frac{2px}{1-\bA(x)}\bigg]=0.
$$
Similarly, we have
$$
G^j(0,y)=\left\{\begin{array}{ll}
\dis\frac{2K_j(y)}{1+\bA(y)} & \mbox{if $j\in\Z^{-*}$,}
\\
0 & \mbox{if $j\in\Z^{+*}$.}
\end{array}\right.
$$
From this and~\refp{generalG}, we immediately derive the function $G(x,y)$.
%
%%%%%%%%%%%%%%
\bth{\label{theorem0}
The generating function $G$ is given by
\beq
G^j(x,y)=\left\{\begin{array}{ll}
\dis\frac{2K_j(x)}{\bA(x)+\bA(y)}& \mbox{if $j\in\Z^{+*}$}
\\[2ex]
\dis\frac{2K_j(y)}{\bA(x)+\bA(y)}& \mbox{if $j\in\Z^{-*}$}
\end{array}\right.
\label{geneGpinned}
\eeq
where $\bA(z)=\sqrt{1-4pqz^2}$.
}
%%%%%%%%%%%%%%
%

%%%%%%%%%%%%%%%%%%%%%%%%%%%%%%%%%%%%%%%%%%%%%%%%%%%%%%%%%%%%%%%%%%%%%%%%%%%%%%%
%%%%%%%%%%%%%%%%%%%%%%%%%         Inverting      %%%%%%%%%%%%%%%%%%%%%%%%%%%%%%
%%%%%%%%%%%%%%%%%%%%%%%%%%%%%%%%%%%%%%%%%%%%%%%%%%%%%%%%%%%%%%%%%%%%%%%%%%%%%%%
\subsection{Distribution of the sojourn time}

We now invert the generating function $G^j$ given by~\refp{geneGpinned}
in order to derive the coefficients $r_{k,n}^j$.
%
%%%%%%%%%%%%%%
\bth{\label{theorem0bis}
The probability $r_{k,n}^j=\P\{T_n=k,S_n=j\}$ admits the following expression:
for $0\le k\le n$ such that $n-k$ is even,
\\
if $j\ge 1$:
\beq
r_{k,n}^j=\left\{\begin{array}{ll}
\dis \!\! 2j \!\!\sum_{i\in\N:\atop {j\le i\le k,\atop i-j\in\cE}}\!\!
\frac{a_{n-i}}{i} \binom{i}{(i+j)/2} p^{(n+j)/2}q^{(n-j)/2}
&\mbox{if $j\le k$ and $k-j$ is even,}
\\
\!\! 0&\mbox{if $j>k$ or $k-j$ is odd;}
\end{array}\right.
\label{law-pinned1}
\eeq
if $j\le -1$:
\begin{align}
r_{k,n}^j=\left\{\begin{array}{ll}
\dis \!\!2|j| \!\!\sum_{i\in\N:\atop {|j|\le i\le k,\atop i+j\in\cE}}
\!\!\frac{a_{n-i}}{i} \binom{i}{(i+j)/2} p^{(n+j)/2}q^{(n-j)/2}
&\mbox{if $j\ge k-n$ and $k-j$ is even,}
\\
\!\! 0&\mbox{if $j<k-n$ or $k-j$ is odd,}
\end{array}\right.
\label{law-pinned2}
\end{align}
where $a_i=\dis\frac{1}{i+2}\binom{i}{i/2}$ for $i\in\cE$.
}
%%%%%%%%%%%%%%
%
\dem
Assume first that $j\ge 1$. We expand $G^j(x,y)$ by using~\refp{1surA+A}:
$$
G^j(x,y)=2\sum_{i=0}^{\infty} \P\{\tau_j=i\}x^i \sum_{l,m\in\cE} a_{l+m}\,x^ly^m
=2\sum_{i,l,m\in\N:\atop l,m\in\cE} a_{l+m}\,\P\{\tau_j=i\}\,x^{i+l}y^m.
$$
By performing the transformations $i+l=k$ and $m=n-k$ in the last sum,
we get
$$
G^j(x,y)=2\sum_{i,k,n\in\N:\atop k-i\in\cE,n-k\in\cE}
a_{n-i}\,\P\{\tau_j=i\}\,x^ky^{n-k}=\sum_{k,n\in\N:\atop n-k\in\cE}
\bigg(\sum_{i\in\N:\atop k-i\in\cE} 2a_{n-i}\,\P\{\tau_j=i\}\bigg)x^ky^{n-k}.
$$
We finally obtain~\refp{law-pinned1} by identifying the coefficients of
the two expansions of $G^j$ and using~\refp{law-tau-a}.
For the case where $j\le -1$, we invoke an argument of duality which
is explained in Remark~\ref{rem-duality} below. The expression of $r_{k,n}^j$
for $j\le -1$ can be deduced from that related to the case where $j\ge 1$ by
interchanging $p$ and $q$, $k$ and $n-k$, $j$ and $-j$, and this
proves~\refp{law-pinned2}.
\fin

%%%%%%%%%%%%%%
\brem{\label{rem-duality}
Let us introduce the dual random walk $(S_n^*)_{n\ge 0}$ with
$S_k^*=\sum_{j=1}^k X_j^*=-S_k$. This sequence is the Bernoulli random walk with
interchanged parameters \mbox{$q=\P\{X_j^*=+1\}$} and $p=\P\{X_j^*=-1\}$.
The corresponding sojourn time is defined as $T_n^*=\sum_{j=1}^n \d_j^*$ with
\begin{align*}
\d_j^*
&=
\left\{\begin{array}{ll}
1 & \mbox{if $(S_j^*>0)$ or $(S_j^*=0$ and $S_{j-1}^*>0)$,}
\\
0 & \mbox{if $(S_j^*<0)$ or $(S_j^*=0$ and $S_{j-1}^*<0)$,}
\end{array}\right.
\\
&=
\left\{\begin{array}{ll}
1 & \mbox{if $(S_j<0)$ or $(S_j=0$ and $S_{j-1}<0)$,}
\\
0 & \mbox{if $(S_j>0)$ or $(S_j=0$ and $S_{j-1}>0)$.}
\end{array}\right.
\end{align*}
We see that $\d_j^*=1-\d_j$ and then $T_n^*=n-T_n$ which implies
$r_{k,n}^j=\P\{\mbox{$T_n^*=n-k,$}\mbox{$S_n^*=-j$}\}.$
As a result, the probability $r_{k,n}^j$ can be deduced from the probability
$r_{n-k,n}^{-j}$ by interchanging $p$ and $q$.
}
%%%%%%%%%%%%%%

%%%%%%%%%%%%%%%%%%%%%%%%%%%%%%%%%%%%%%%%%%%%%%%%%%%%%%%%%%%%%%%%%%%%%%%%%%%%%%%
%%%%%%%%%%%%%%%%%%%%%%%%%     Random walk      %%%%%%%%%%%%%%%%%%%%%%%%%%%%%%%%
%%%%%%%%%%%%%%%%%%%%%%%%%%%%%%%%%%%%%%%%%%%%%%%%%%%%%%%%%%%%%%%%%%%%%%%%%%%%%%%
\section{The case $F=\Z$ (random walk without conditioning)}

In this part, we study the sojourn time without conditioning the extremity
of the random walk. This corresponds to the case $F=\Z$. We set for simplifying
the notations \mbox{$r_{k,n}^{\Z}=r_{k,n}$} and $G^{\Z}(x,y)=G(x,y)$.
A possible expression for the $r_{k,n}$'s can be obtained from
Theorem~\ref{theorem0bis} by summing the $r_{k,n}^j$, $j\in\Z$. Actually,
$r_{k,n}^j=0$ for $|j|>n$; so,
$$
r_{k,n}=\sum_{j=-n}^{n} r_{k,n}^j.
$$
We propose another representation which can be deduced from the generating function $G$.

%%%%%%%%%%%%%%%%%%%%%%%%%%%%%%%%%%%%%%%%%%%%%%%%%%%%%%%%%%%%%%%%%%%%%%%%%%%%%%%
%%%%%%%%%%%%%%%%%%%%%%%%%  Generating function   %%%%%%%%%%%%%%%%%%%%%%%%%%%%%%
%%%%%%%%%%%%%%%%%%%%%%%%%%%%%%%%%%%%%%%%%%%%%%%%%%%%%%%%%%%%%%%%%%%%%%%%%%%%%%%
\subsection{Generating function}

In order to derive the generating function $G(x,y)$, in view of~\refp{generalG},
we need to evaluate the functions $H_i^\Z(z)$, $G(x,0)$ and $G(0,y)$.
On one hand, by~\refp{geneH},
\begin{align*}
H_i^\Z(z) &= \sum_{j\in \Z-i} H_j(z)=\sum_{j\in\Z} H_j(z)
\\
&=
\frac{1}{\bA(z)}\bigg[\sum_{j=1}^{\infty} \bigg(\frac{1-\bA(z)}{2pz}\bigg)^j+
\sum_{j=1}^{\infty} \bigg(\frac{1-\bA(z)}{2qz}\bigg)^j+1\bigg]
\\
&=
\frac{1}{\bA(z)}\bigg[\frac{1-\bA(z)}{2pz}\,\frac{1}{1-\frac{1-\bA(z)}{2pz}}
+\frac{1-\bA(z)}{2qz}\,\frac{1}{1-\frac{1-\bA(z)}{2qz}}+1\bigg]
\\
&=
\frac{1}{\bA(z)}\bigg[\frac{1-\bA(z)}{2pz-1+\bA(z)}
+\frac{1-\bA(z)}{2qz-1+\bA(z)}+1\bigg].
\end{align*}
We have
\begin{align}
[2pz-1+\bA(z)][2qz-1+\bA(z)]
&=
4pqz^2-2z+1+2(z-1)\bA(z)+\bA(z)^2
\nonumber\\
&=
2(1-z)[1-\bA(z)].
\label{identityA}
\end{align}
Then
$$
H_i^\Z(z) =\frac{1}{\bA(z)}\bigg[\frac{2qz-1+\bA(z)}{2(1-z)}
+\frac{2pz-1+\bA(z)}{2(1-z)}+1\bigg].
$$
The term within the brackets equals $\bA(z)/(1-z)$ and thus, for any $i\in\Z$,
$$
H_i^\Z(z)=\frac{1}{1-z}.
$$
Notice that the value of $H_i^\Z(z)$ may be deduced through a shorter way from
its own definition:
$$
H_i^\Z(z)=\sum_{k=0}^{\infty} \P_i\{S_k\in\Z\} z^k
=\sum_{k=0}^{\infty} z^k=\frac{1}{1-z}.
$$

On the other hand, by~\refp{GF(x,0)},
$$
G(x,0)=H_0^\Z(x)-\frac{1-\bA(x)}{2px}\,H_{-1}^\Z(x)
=\frac{1}{1-x}\bigg[1-\frac{1-\bA(x)}{2px}\bigg]
=\frac{2px-1+\bA(x)}{2px(1-x)}.
$$
Similarly,
$$
G(0,y)=\frac{2qy-1+\bA(y)}{2qy(1-y)}.
$$
We can then write out $G(x,y)$.
%
%%%%%%%%%%%%%%
\bth{\label{theorem1}
The generating function $G$ is given by
\beq
G(x,y)=\frac{1}{\bA(x)+\bA(y)}\bigg[\frac{2p-1+\bA(x)}{1-x}
+\frac{2q-1+\bA(y)}{1-y}\bigg]
\label{geneGfree}
\eeq
where $\bA(z)=\sqrt{1-4pqz^2}$.
}
%%%%%%%%%%%%%%
%
\dem
We have by~\refp{generalG}
$$
G(x,y)=\frac{[1+\bA(x)]G(x,0)+[1+\bA(y)]G(0,y)-2}{\bA(x)+\bA(y)}.
$$
Let us compute the terms $[1+\bA(x)]G(x,0)$ and $[1+\bA(y)]G(0,y)$. We have
$$
[1+\bA(x)][2px-1+\bA(x)]=2px-1+2px\bA(x)+\bA(x)^2=2px[1-2qx+\bA(x)]
$$
and then
$$
[1+\bA(x)]G(x,0)=\frac{1-2qx+\bA(x)}{1-x}.
$$
Similarly,
$$
[1+\bA(y)]G(0,y)=\frac{1-2py+\bA(y)}{1-y}.
$$
Therefore, we get
\begin{align*}
\lqn{[1+\bA(x)]G(x,0)+[1+\bA(y)]G(0,y)-2}
\\[-3ex]
&=
\bigg[\frac{1-2qx}{1-x}+\frac{1-2py}{1-y}-2\bigg]+\frac{\bA(x)}{1-x}
+\frac{\bA(y)}{1-y}
\\
&=
\bigg[\bigg(\frac{1-2qx}{1-x}-2q\bigg)
+\bigg(\frac{1-2py}{1-y}-2p\bigg)\bigg]+\frac{\bA(x)}{1-x}+\frac{\bA(y)}{1-y}
\\
&=
\frac{2p-1+\bA(x)}{1-x}+\frac{2q-1+\bA(y)}{1-y}
\end{align*}
from which we obtain~\refp{geneGfree}.
\fin

An interesting consequence of Theorem~\ref{theorem1} concerns the ``partial''
generating function $\tG$ of the probabilities $r_{k,n}$ limited to the even
indices $k$ and $n$:
$$
\tG(x,y)=\sum_{k,n\in\cE} r_{k,n} x^ky^{n-k}.
$$
For this function, we have the simple result below.
%
%%%%%%%%%%%%%%
\bco{\label{corollary1}
The generating function $\tG$ is given by
\beq
\tG(x,y)=\frac{4pq}{[1-2p+\bA(x)][1-2q+\bA(y)]}
\label{geneGfree-even}
\eeq
where $\bA(z)=\sqrt{1-4pqz^2}$. In particular,
\beq
\tG(x,y)=\tG(x,0)\tG(0,y).
\label{geneGfree-even-product}
\eeq
}
%%%%%%%%%%%%%%
%
\dem
We first try to relate $\tG$ to $G$. For this, we observe that
$$
\frac12[G(x,y)+G(-x,-y)]
=\sum_{k,n\in\N:\atop k\le n} \frac{1+(-1)^n}{2}\,r_{k,n} x^ky^{n-k}
=\sum_{k\in\N,n\in\cE:\atop k\le n} r_{k,n} x^ky^{n-k}.
$$
We know that, for all even integer $n$, if $k$ is odd, then $r_{k,n}=0$. We
therefore have checked that
$$
\tG(x,y)=\frac12[G(x,y)+G(-x,-y)].
$$
Thus,
$$
\tG(x,y)=\frac{1}{\bA(x)+\bA(y)}\bigg[\frac{2p-1+\bA(x)}{1-x^2}
+\frac{2q-1+\bA(y)}{1-y^2}\bigg].
$$
We have now
$$
2p-1+\bA(x)=\frac{(2p-1)^2-\bA(x)^2}{2p-1-\bA(x)}=4pq\,\frac{1-x^2}{1-2p+\bA(x)}
$$
and similarly,
$$
2q-1+\bA(y)=4pq\,\frac{1-y^2}{1-2q+\bA(y)}.
$$
This gives
$$
\tG(x,y)=\frac{4pq}{\bA(x)+\bA(y)}\bigg[\frac{1}{1-2p+\bA(x)}
+\frac{1}{1-2q+\bA(y)}\bigg]
$$
which in turn, with $[1-2p+\bA(x)]+[1-2q+\bA(y)]=\bA(x)+\bA(y),$
yields~\refp{geneGfree-even}. This formula supplies in particular
$$
\tG(x,0)=\frac{2q}{1-2p+\bA(x)} \quad\mbox{and}\quad \tG(0,y)=\frac{2p}{1-2q+\bA(y)}
$$
and we immediately obtain~\refp{geneGfree-even-product}.
\fin

%%%%%%%%%%%%%%%%%%%%%%%%%%%%%%%%%%%%%%%%%%%%%%%%%%%%%%%%%%%%%%%%%%%%%%%%%%%%%%%
%%%%%%%%%%%%%%%%%%%%%%%%%         Inverting      %%%%%%%%%%%%%%%%%%%%%%%%%%%%%%
%%%%%%%%%%%%%%%%%%%%%%%%%%%%%%%%%%%%%%%%%%%%%%%%%%%%%%%%%%%%%%%%%%%%%%%%%%%%%%%
\subsection{Distribution of the sojourn time}\label{section-free}

In this part, we invert the generating function $G$ given by~\refp{geneGfree}
in order to derive the coefficients $r_{k,n}$.
%
%%%%%%%%%%%%%%
\bth{\label{theorem1bis}
The probability $r_{k,n}=\P\{T_n=k\}$ admits the following expression:
for $0\le k\le n$,
\begin{align}
r_{k,n}
&=
\ind_{\cE}(n-k)\bigg[2p\sum_{i\in\cE:\atop n-k\le i\le n} a_i (pq)^{i/2}
-4\sum_{i\in\cE:\atop n-k\le i\le n-2}
\bigg(\sum_{j\in \cE:\atop j\le i+k-n} a_ja_{i-j}\bigg) (pq)^{i/2+1}\bigg]
\nonumber\\
&\hphantom{=\,}
+\ind_{\cE}(k)\bigg[2q\sum_{i\in\cE:\atop k\le i\le n} a_i (pq)^{i/2}
-4\sum_{i\in\cE:\atop k\le i\le n-2}
\bigg(\sum_{j\in \cE:\atop j\le i-k} a_ja_{i-j}\bigg) (pq)^{i/2+1}\bigg]
\label{lawfree}
\end{align}
where $a_i=\dis\frac{1}{i+2}\binom{i}{i/2}$ for $i\in\cE$.
More specifically,
\bitem
\item
For odd $n$:
\beq
r_{k,n}=\left\{\begin{array}{l}
\dis 2p\sum_{i\in\cE:\atop n-k\le i\le n-1} a_i (pq)^{i/2}
-4\sum_{i\in\cE:\atop n-k\le i\le n-3}
\bigg(\sum_{j\in \cE:\atop j\le i+k-n} a_ja_{i-j}\bigg) (pq)^{i/2+1}
\\
\mbox{if $k$ is odd,}
\\[1ex]
\dis 2q\sum_{i\in\cE:\atop k\le i\le n-1} a_i (pq)^{i/2}
-4\sum_{i\in\cE:\atop k\le i\le n-3}
\bigg(\sum_{j\in \cE:\atop j\le i-k} a_ja_{i-j}\bigg) (pq)^{i/2+1}
\\
\mbox{if $k$ is even;}
\end{array}\right.
\label{lawfree-odd-case}
\eeq
\item
For even $n$:
\beq
r_{k,n}=\left\{\begin{array}{l}
\dis 2p\sum_{i\in\cE:\atop n-k\le i\le n} a_i (pq)^{i/2}
-4\sum_{i\in\cE:\atop n-k\le i\le n-2}
\bigg(\sum_{j\in \cE:\atop j\le i+k-n} a_ja_{i-j}\bigg) (pq)^{i/2+1}
\\[4ex]
\dis +2q\sum_{i\in\cE:\atop k\le i\le n} a_i (pq)^{i/2}
-4\sum_{i\in\cE:\atop k\le i\le n-2}
\bigg(\sum_{j\in \cE:\atop j\le i-k} a_ja_{i-j}\bigg) (pq)^{i/2+1}
\\
\mbox{if $k$ is even,}
\\[1ex]
0\mbox{ if $k$ is odd.}
\end{array}\right.
\label{lawfree-even-case}
\eeq
\eitem
}
%%%%%%%%%%%%%%
%
\dem
Rewrite~\refp{geneGfree} as
\begin{align}
G(x,y)
&=
\frac{2p}{1-x}\,\frac{1}{\bA(x)+\bA(y)}
-\frac{1-\bA(x)}{1-x}\,\frac{1}{\bA(x)+\bA(y)}
\nonumber\\
&\hphantom{=\,}
+\frac{2q}{1-y}\,\frac{1}{\bA(x)+\bA(y)}
-\frac{1-\bA(y)}{1-y}\,\frac{1}{\bA(x)+\bA(y)}.
\label{G4terms}
\end{align}

We expand the first term of~\refp{G4terms}. By~\refp{elem} and~\refp{1surA+A},
\begin{align}
\frac{2p}{1-x}\,\frac{1}{\bA(x)+\bA(y)}
&=
2p\bigg(\sum_{m\in\N} x^m\bigg)
\bigg(\sum_{j,l\in \cE} a_{j+l} (pq)^{(j+l)/2}x^jy^l\bigg)
\nonumber\\
&=
2p\sum_{j,l\in\cE,m\in\N} a_{j+l}(pq)^{(j+l)/2}x^{j+m}y^l
\nonumber\\
&=
2p\sum_{k,n\in\N} \bigg[\sum_{j,l\in \cE,m\in\N:
\atop j+m=k,l=n-k} a_{j+l}(pq)^{(j+l)/2}\bigg] x^ky^{n-k}
\label{first-term-intermediate1}
\\
&=
2p\sum_{k,n\in\N:\atop k\le n,n-k\in\cE} \bigg[\sum_{m\in\N:\atop m \le k,n-m\in\cE}
a_{n-m} (pq)^{(n-m)/2}\bigg] x^ky^{n-k}
\label{first-term-intermediate2}
\\
&=
2p\sum_{k,n\in\N:\atop k\le n,} \ind_{\cE}(n-k)
\bigg[\sum_{i\in\cE:\atop n-k\le i\le n} a_i (pq)^{i/2}\bigg] x^ky^{n-k}.
\label{first-termG}
\end{align}
Let us explain certain transformations made in the above calculations:
\begin{enumerate}
\item
In~\refp{first-term-intermediate1}, we have introduced the indices
$k=j+m$ and $n=j+l+m$. Then $j=k-m$ and $l=n-k$,
the conditions $j,l\in\cE$ are equivalent to
$k-m\ge0,n-k\ge0,k-m\in\cE,$ $n-k\in\cE$,
or $m\le k\le n, n-k\in\cE, n-m\in\cE$. This leads
to~\refp{first-term-intermediate2};
\item
From~\refp{first-term-intermediate2} to~\refp{first-termG}, we have put $i=n-m$.
\end{enumerate}

We expand the second term of~\refp{G4terms}. By~\refp{elem} and~\refp{1surA+A},
\begin{align}
\lqn{\frac{1-\bA(x)}{1-x}\,\frac{1}{\bA(x)+\bA(y)}}
\nonumber\\[-3ex]
&=
4\bigg(\sum_{m\in\N} x^m\bigg) \bigg(\sum_{i\in \cE^*} a_{i-2} (pq)^{i/2}x^i\bigg)
\bigg(\sum_{j,l\in \cE} a_{j+l} (pq)^{(j+l)/2}x^jy^l\bigg)
\nonumber\\
&=
4\sum_{i\in \cE^*\!,j,l\in\cE,m\in\N} a_{i-2}a_{j+l} (pq)^{(i+j+l)/2}x^{i+j+m}y^l
\nonumber\\
&=
4\sum_{k,n\in\N} \bigg[\sum_{i\in \cE^*\!,j,l\in\cE,m\in\N:
\atop i+j+m=k,l=n-k} a_{i-2}a_{j+l} (pq)^{(i+j+l)/2}\bigg] x^ky^{n-k}
\label{second-term-intermediate1}
\\
\lqn{}
&=
4\sum_{k,n\in\N:\atop k\le n,n-k\in\cE} \bigg[\sum_{i\in \cE^*\!,m\in\N:
\atop i+m\le k,i+n-m\in\cE} a_{i-2}a_{n-m-i} (pq)^{(n-m)/2}\bigg] x^ky^{n-k}
\label{second-term-intermediate2}
\\
&=
4\sum_{k,n\in\N:\atop k\le n,n-k\in\cE} \bigg[\sum_{m\in\N:\atop n-m\in\cE}
\bigg(\sum_{i\in \cE^*\!:\atop i\le k-m} a_{i-2}a_{n-m-i}\bigg)
(pq)^{(n-m)/2}\bigg] x^ky^{n-k}
\label{second-term-intermediate3}
\\
&=
4\sum_{k,n\in\N:\atop k\le n,n-k\in\cE} \bigg[\sum_{m\in\N:\atop n-m\in\cE}
\bigg(\sum_{j\in \cE:\atop j\le k-m-2} a_j a_{n-m-j-2}\bigg)
(pq)^{(n-m)/2}\bigg] x^ky^{n-k}
\label{second-term-intermediate4}
\\
&=
4\sum_{k,n\in\N:\atop k\le n,n-k\in\cE}
\bigg[\sum_{m\in\N\setminus\{0,1\}:\atop n-m\in\cE}
\bigg(\sum_{j\in \cE:\atop j\le k-m} a_ja_{n-m-j}\bigg)
(pq)^{(n-m)/2+1}\bigg] x^ky^{n-k}
\label{second-term-intermediate5}
\\
&=
4\sum_{k,n\in\N:\atop k\le n} \ind_{\cE}(n-k)\bigg[\sum_{i\in\cE:\atop i\le n-2}
\bigg(\sum_{j\in \cE:\atop j\le i+k-n} a_ja_{i-j}\bigg)
(pq)^{i/2+1}\bigg] x^ky^{n-k}.
\label{second-termG}
\end{align}
Let us explain the transformations used in the above computations:
\begin{enumerate}
\item
In~\refp{second-term-intermediate1}, we have introduced the indices
$k=i+j+m$ and $n=i+j+l+m$. Then $j=-i+k-m,l=n-k$ and $i\in\cE$,
the conditions $j,l\in\cE$ are equivalent to
$k-i-m\ge0,n-k\ge0,k-i-m\in\cE,n-k\in\cE$,
or $i+m\le k\le n, n-k\in\cE, n-m\in\cE$. This leads
to~\refp{second-term-intermediate2};
\item
From~\refp{second-term-intermediate3} to~\refp{second-term-intermediate4},
we have put $i=j+2$;
\item
From~\refp{second-term-intermediate4} to~\refp{second-term-intermediate5},
we have changed $m$ into $m-2$;
\item
From~\refp{second-term-intermediate5} to~\refp{second-termG}, we have put $i=n-m$.
\end{enumerate}

By invoking the argument of duality which is explained in Remark~\ref{rem-duality},
the two last terms of~\refp{G4terms} can be deduced from the two first ones
by interchanging $p$ and $q$, and $x$ and $y$ (that is $k$ and $n-k$). This yields
\begin{align}
\frac{2q}{1-y}\,\frac{1}{\bA(x)+\bA(y)}
&=
2q\sum_{k,n\in\N:\atop k\le n} \ind_{\cE}(k)\bigg[\sum_{i\in\cE:\atop k\le i\le n}
a_i (pq)^{i/2}\bigg] x^ky^{n-k},
\label{third-termG}
\\
\frac{1-\bA(y)}{1-y}\,\frac{1}{\bA(x)+\bA(y)}
&=
4\sum_{k,n\in\N:\atop k\le n}  \ind_{\cE}(k)\bigg[\sum_{i\in\cE:\atop i\le n-2}
\bigg(\sum_{j\in \cE:\atop j\le i-k} a_ja_{i-j}\bigg) (pq)^{i/2+1}\bigg] x^ky^{n-k}.
\label{fourth-termG}
\end{align}
By identifying the coefficients of the series lying in the definition of $G$
and in \refp{first-termG}, \refp{second-termG}, \refp{third-termG} and
\refp{fourth-termG}, we extract~\refp{lawfree}.
\fin

%
%%%%%%%%%%%%%%
\bco{
The probabilities
\begin{align*}
r_{0,n}&=\P\{T_n=0\}=\P\{S_1\le 0,\dots,S_n\le 0\}=\P\{\tau_1>n\},
\\
r_{n,n}&=\P\{T_n=n\}=\P\{S_1\ge 0,\dots,S_n\ge 0\}=\P\{\tau_{-1}>n\}
\end{align*}
are given by
\beq
r_{0,n}=1-2p\sum_{i\in\cE:\atop i\le n-1} a_i (pq)^{i/2}\quad\mbox{and}\quad
r_{n,n}=1-2q\sum_{i\in\cE:\atop i\le n-1} a_i (pq)^{i/2}.
\label{lawfree-part-cases}
\eeq
The probabilities
\begin{align*}
r_{1,n}&=\P\{T_n=1\}=\P\{S_1\le 0,\dots,S_{n-1}\le 0,S_n>0\}=\P\{\tau_1=n\},
\\
r_{n-1,n}&=\P\{T_n=n-1\}=\P\{S_1\ge 0,\dots,S_{n-1}\ge 0,S_n<0\}=\P\{\tau_{-1}=n\}
\end{align*}
are given by
$r_{1,n}=r_{n-1,n}=0$ if $n$ is even and, if $n$ is odd,
\beq
r_{1,n}= 2 a_{n-1}\, p^{(n+1)/2}q^{(n-1)/2} \quad\mbox{and}\quad
r_{n-1,n}=2 a_{n-1}\, p^{(n-1)/2}q^{(n+1)/2}.
\label{lawfree-part-casesbis}
\eeq
}
%%%%%%%%%%%%%%
%
We propose three proofs of~\refp{lawfree-part-cases}.

\demun
Formulas~\refp{lawfree-part-casesbis} can be easily deduced
from~\refp{lawfree-odd-case} and~\refp{lawfree-even-case}.
Now, we deduce~\refp{lawfree-part-cases} from Theorem~\ref{theorem1bis}.
We only consider the probability $r_{n,n}$, the computations related to
$r_{0,n}$ being quite analogous.

If $n$ is odd, \refp{lawfree-odd-case} gives
\begin{align*}
r_{n,n}
&=
2p\sum_{i\in\cE:\atop i\le n-1} a_i (pq)^{i/2}
-4\sum_{i\in\cE:\atop i\le n-3}
\bigg(\sum_{j\in \cE:\atop j\le i} a_ja_{i-j}\bigg) (pq)^{i/2+1}
\\
&=
(2-2q)\sum_{i\in\cE:\atop i\le n-1} a_i (pq)^{i/2}
-2\sum_{i\in\cE:\atop i\le n-3} a_{i+2} (pq)^{i/2+1}
\\
&=
(2-2q)\sum_{i\in\cE:\atop i\le n-1} a_i (pq)^{i/2}
-2\sum_{i\in\cE^*:\atop i\le n-1} a_i (pq)^{i/2}
\\
&=
1-2q\sum_{i\in\cE:\atop i\le n-1} a_i (pq)^{i/2}.
\end{align*}

If $n$ is even, \refp{lawfree-even-case} gives
\begin{align*}
r_{n,n}
&=
2p\sum_{i\in\cE:\atop i\le n} a_i (pq)^{i/2}
+2q a_n (pq)^{i/2}-4\sum_{i\in\cE:\atop i\le n-2}
\bigg(\sum_{j\in \cE:\atop j\le i} a_ja_{i-j}\bigg) (pq)^{i/2+1}
\\
&=
(2-2q)\sum_{i\in\cE:\atop i\le n} a_i (pq)^{i/2}
+2q a_n (pq)^{i/2}-2\sum_{i\in\cE:\atop i\le n-2} a_{i+2} (pq)^{i/2+1}
\\
&=
2\sum_{i\in\cE:\atop i\le n} a_i (pq)^{i/2}
-2q\sum_{i\in\cE:\atop i\le n-2} a_i (pq)^{i/2}
-2\sum_{i\in\cE^*:\atop i\le n} a_i (pq)^{i/2}
\\
&=
1-2q\sum_{i\in\cE:\atop i\le n-2} a_i (pq)^{i/2}
=1-2q\sum_{i\in\cE:\atop i\le n-1} a_i (pq)^{i/2}.
\end{align*}
\fin

\demdeux
Another proof of~\refp{lawfree-part-cases} consists of considering the
generating function of the $r_{n,n}$'s, $n\in\N$, which is nothing
but $G(x,0)$ as observed in~\refp{partial}. Indeed, we have by~\refp{elem}
\begin{align*}
G(x,0)
&=
\frac{2px-1+\bA(x)}{2px(1-x)}=\frac{1}{1-x}-\frac{1-\bA(x)}{2px(1-x)}.
\\
&=
\sum_{m\in\N} x^m-2q\bigg(\sum_{m\in\N} x^m\bigg)
\bigg(\sum_{i\in \cE} a_i (pq)^{i/2}x^{i+1}\bigg)
\\
\lqn{}
&=
\sum_{m\in\N} x^m-2q\sum_{i\in \cE,m\in\N} a_i (pq)^{i/2}x^{i+m+1}
\\
&=
\sum_{n\in\N} x^n-2q\sum_{n\in\N^*} \bigg(\sum_{i\in \cE:\atop i\le n-1}
a_i (pq)^{i/2}\bigg) x^n
\\
&=
\sum_{n\in\N} \bigg(1-2q\sum_{i\in \cE:\atop i\le n-1} a_i (pq)^{i/2}\bigg) x^n.
\end{align*}
From this and~\refp{partial}, we obtain the expression~\refp{lawfree-part-cases}
related to $r_{n,n}$ by identification.
\fin

\demtrois
A direct proof for deriving~\refp{lawfree-part-cases} consists of writing,
by using~\refp{law-tau1}, that
\begin{align*}
r_{n,n}
&=
\P\{\tau_{-1}>n\}=1-\P\{\tau_{-1}\le n\}
=1-\sum_{k\in\cO:\atop k\le n} \P\{\tau_{-1}=k\}
\\
&=
1-\frac{1}{2p}\sum_{k\in\cO:\atop k\le n} \pi_{k+1}
=1-\frac{2}{p}\sum_{k\in\cO:\atop k\le n} a_{k-1} (pq)^{(k+1)/2}
=1-2q\sum_{k\in\cE:\atop k\le n-1} a_k (pq)^{k/2}.
\end{align*}
\fin

%
%%%%%%%%%%%%%%
\brem{
Let us compute the limit of $r_{0,n}$ as $n$ tends to $\infty$:
\refp{lawfree-part-cases} gives
$$
\lim_{n\to \infty} r_{0,n}=1-2p\sum_{i\in\cE} a_i (pq)^{i/2}.
$$
By~\refp{elem}, one has
\beq
\sum_{i\in\cE} a_i (pq)^{i/2}=\frac{1-A(\sqrt{pq})}{4pq}=\frac{1-|p-q|}{4pq}
=\left\{\begin{array}{ll}
\dis\frac{1}{2p} & \mbox{if $p\ge 1/2$}
\\[3ex]
\dis\frac{1}{2q} & \mbox{if $p\le 1/2$}
\end{array}\right.
\label{elembis}
\eeq
which yields
\beq
\lim_{n\to \infty} r_{0,n}=\left\{\begin{array}{ll}
0 & \mbox{if $p\ge 1/2$,}
\\[1ex]
\dis 2-\frac 1q & \mbox{if $p\le 1/2$.}
\end{array}\right.
\label{limit-r0n}
\eeq
This simple result is in good accordance with the fact that
$$
\lim_{n\to \infty} r_{0,n}=\P\{T_{\infty}=0\}=\P\{\tau_1=\infty\}.
$$
Analogously,
$$
\lim_{n\to \infty} r_{n,n}=\P\{\tau_{-1}=\infty\}
=\left\{\begin{array}{ll}
0 & \mbox{if $p\le 1/2$,}
\\[1ex]
\dis 2-\frac 1p & \mbox{if $p\ge 1/2$.}
\end{array}\right.
$$
}
%%%%%%%%%%%%%%
%

%
%%%%%%%%%%%%%%
\bco{
For even $n,k$ such that $0\le k\le n$, the following relationship holds:
\beq
r_{k,n}=r_{k,k}r_{0,n-k}.
\label{law-product}
\eeq
}
%%%%%%%%%%%%%%
%
We propose three proofs of~\refp{law-product}. The first one is based on the
explicit representation~\refp{lawfree-part-cases} of the probabilities $r_{k,k}$
and $r_{0,n-k}$ and the representation~\refp{lawfree-even-case}
of $r_{k,n}$ (so, this is a consequence of Theorem~\ref{theorem1bis}),
the second one is done by induction and the third one relies on the specific
generating function $\tG$ which has been introduced below Theorem~\ref{theorem1}.

\demun
We deduce~\refp{law-product} from~\refp{lawfree-part-cases}.
Let us develop the product $r_{k,k}r_{0,n-k}$. The computations are
rather technical:
\begin{align}
r_{k,k}r_{0,n-k}
&=
\bigg[1-2p\sum_{i\in\cE:\atop i\le n-k-2} a_i (pq)^{i/2}\bigg]
\bigg[1-2q\sum_{i\in\cE:\atop i\le k-2} a_i (pq)^{i/2}\bigg]
\nonumber\\
&=
1-2p\sum_{i\in\cE:\atop i\le n-k-2} a_i (pq)^{i/2}
-2q\sum_{i\in\cE:\atop i\le k-2} a_i (pq)^{i/2}
\nonumber\\
&\hphantom{=\,}
+4pq \sum_{l,m\in\cE:\atop l\le k-2,m\le n-k-2} a_la_m (pq)^{(l+m)/2}.
\label{product-intermediate1}
\end{align}
We write the two first sums of the last displayed equation as
\begin{align}
\sum_{i\in\cE:\atop i\le n-k-2} a_i (pq)^{i/2}
&=
\sum_{i\in\cE:\atop i\le n} a_i (pq)^{i/2}
-\sum_{i\in\cE:\atop n-k \le i\le n} a_i (pq)^{i/2},
\nonumber\\[-1ex]
\label{product-intermediate2}
\\[-1ex]
\sum_{i\in\cE:\atop i\le k-2} a_i (pq)^{i/2}
&=
\sum_{i\in\cE:\atop i\le n} a_i (pq)^{i/2}
-\sum_{i\in\cE:\atop k \le i\le n} a_i (pq)^{i/2}
\nonumber
\end{align}
and the third one as
\begin{align}
\sum_{l,m\in\cE:\atop l\le k-2,m\le n-k-2} a_la_m (pq)^{(l+m)/2}
&=
\sum_{i\in\cE:\atop i\le n-4}\bigg(\sum_{j\in\cE:\atop i+k+2-n\le j\le (k-2)\wedge i}
a_ja_{i-j}\bigg) (pq)^{i/2}
\nonumber\\
&=
\sum_{i\in\cE:\atop i\le k-2}\bigg(\sum_{j\in\cE:\atop j\le i}
a_ja_{i-j}\bigg) (pq)^{i/2}
+\!\!\!\!\sum_{i\in\cE:\atop k\le i\le n-k-2}\!\!\!\!\bigg(\sum_{j\in\cE:\atop j\le k-2}
a_ja_{i-j}\bigg) (pq)^{i/2}
\nonumber\\
&\hphantom{=\,}
+\sum_{i\in\cE:\atop n-k\le i\le n-4}\bigg(\sum_{j\in\cE:\atop i+k+2-n\le j\le k-2}
a_ja_{i-j}\bigg) (pq)^{i/2}.
\label{product-intermediate3}
\end{align}
In the first equality of~\refp{product-intermediate2}, we changed the indices
$l,m$ into $i,j$ according to the rules $i=l+m$ and $j=l$. The conditions
$l\le k-2$ and $m\le n-k-2$ yield $i\le n-4$ and $i+k+2-n\le j\le i$
(since $j=i-m$). In the second equality of~\refp{product-intermediate2},
we split the sum into three terms: for $i\le k-2$, the condition
$i+k+2-n\le j\le (k-2)\wedge i$ is equivalent to $j\le i$,
for $k\le i\le n-k-2$, this condition is equivalent to $j\le k-2$ and
for $n-k\le i\le n-4$, this condition is equivalent to $i+k+2-n\le j\le k-2$.

Therefore, putting~\refp{product-intermediate2} and~\refp{product-intermediate3}
into~\refp{product-intermediate1},
\begin{align}
\lqn{r_{k,k}r_{0,n-k}}
\nonumber\\[-3ex]
&=
2a_0-2p\bigg[\sum_{i\in\cE:\atop i\le n} a_i (pq)^{i/2}-
\sum_{i\in\cE:\atop n-k \le i\le n} a_i (pq)^{i/2}\bigg]
-2q\bigg[\sum_{i\in\cE:\atop i\le n} a_i (pq)^{i/2}-
\sum_{i\in\cE:\atop k \le i\le n} a_i (pq)^{i/2}\bigg]
\nonumber\\
\lqn{}&\hphantom{=\,}
+\,4pq\bigg[\sum_{i\in\cE:\atop i\le k-2}\bigg(\sum_{j\in\cE:\atop j\le i}
a_ja_{i-j}\bigg) (pq)^{i/2}
+\sum_{i\in\cE:\atop k\le i\le n-k-2}\bigg(\sum_{j\in\cE:\atop j\le k-2}
a_ja_{i-j}\bigg) (pq)^{i/2}
\nonumber\\
&\hphantom{=\,}
+\sum_{i\in\cE:\atop n-k\le i\le n-4}\bigg(\sum_{j\in\cE:\atop i+k+2-n\le j\le k-2}
a_ja_{i-j}\bigg) (pq)^{i/2}\bigg]
\nonumber\\
&=
2p\sum_{i\in\cE:\atop n-k \le i\le n} a_i (pq)^{i/2}
+2q\sum_{i\in\cE:\atop k \le i\le n} a_i (pq)^{i/2}
-2\sum_{i\in\cE^*:\atop i\le n} a_i (pq)^{i/2}
+4pq\bigg[\dots\bigg].
\nonumber\\
&=
2p\sum_{i\in\cE:\atop n-k \le i\le n} a_i (pq)^{i/2}
+2q\sum_{i\in\cE:\atop k \le i\le n} a_i (pq)^{i/2}
+4pq\bigg[\dots-\frac12\sum_{i\in\cE^*:\atop i\le n} a_i (pq)^{i/2-1}\bigg].
\label{intermediate-product}
\end{align}
In the second equality of~\refp{intermediate-product}, we have simplified
$a_0-p\sum_{i\in\cE:\atop i\le n} a_i (pq)^{i/2}-
q\sum_{i\in\cE:\atop i\le n} a_i (pq)^{i/2}$
into
$-\sum_{i\in\cE^*:\atop i\le n} a_i (pq)^{i/2}$,
and the dots  ``$\dots$'' stand for the whole term lying within brackets
which follows the factor $4pq$ in the first equality.

Using~\refp{convol}, we rewrite the last sum lying in~\refp{intermediate-product} as
$$
\frac12\sum_{i\in\cE^*:\atop i\le n} a_i (pq)^{i/2-1}
=\frac12\sum_{i\in\cE:\atop i\le n-2} a_{i+2} (pq)^{i/2}
=\sum_{i\in\cE:\atop i\le n-2} \bigg(\sum_{j\in\cE:\atop j\le i}
a_ja_{i-j}\bigg) (pq)^{i/2}
$$
and then \refp{intermediate-product}~becomes
\begin{align*}
r_{k,k}r_{0,n-k}
&=
2p\sum_{i\in\cE:\atop n-k \le i\le n} a_i (pq)^{i/2}
+2q\sum_{i\in\cE:\atop k \le i\le n} a_i (pq)^{i/2}
\\
&\hphantom{=\,}
+4pq\bigg[\sum_{i\in\cE:\atop i\le k-2}\bigg(\sum_{j\in\cE:\atop j\le i}
a_ja_{i-j}\bigg) (pq)^{i/2}
+\sum_{i\in\cE:\atop k\le i\le n-k-2}\bigg(\sum_{j\in\cE:\atop j\le k-2}
a_ja_{i-j}\bigg) (pq)^{i/2}
\\
&\hphantom{=\,}
+\sum_{i\in\cE:\atop n-k\le i\le n-4}\bigg(\sum_{j\in\cE:\atop i+k+2-n\le j\le k-2}
a_ja_{i-j}\bigg) (pq)^{i/2}
-\sum_{i\in\cE:\atop i\le n-2} \bigg(\sum_{j\in\cE:\atop j\le i}a_ja_{i-j}\bigg)
(pq)^{i/2}\bigg].
\end{align*}
We compute the term within brackets. Since the computations
are technical, we first write down progressively all the equalities and we next
explain them. The aforementioned term writes, by omitting the repetitive
quantity $a_ja_{i-j}(pq)^{i/2}$ for lightening the proof,
\begin{align}
&
\sum_{i\in\cE:\atop i\le k-2}\sum_{j\in\cE:\atop j\le i}
+\sum_{i\in\cE:\atop k\le i\le n-k-2}\sum_{j\in\cE:\atop j\le k-2}
+\sum_{i\in\cE:\atop n-k\le i\le n-4}\sum_{j\in\cE:\atop i+k-n+2\le j\le k-2}
-\sum_{i\in\cE:\atop i\le n-2}\sum_{j\in\cE:\atop j\le i}
\label{computation1}\\
&=
-\sum_{i\in\cE:\atop k\le i\le n-2}\sum_{j\in\cE:\atop j\le i}
+\sum_{i\in\cE:\atop k\le i\le n-k-2}\sum_{j\in\cE:\atop j\le k-2}
+\sum_{i\in\cE:\atop n-k\le i\le n-4}\sum_{j\in\cE:\atop i-k+2\le j\le n-k-2}
\label{computation2}\\
&=
-\sum_{i\in\cE:\atop k\le i\le n-2}\sum_{j\in\cE:\atop j\le k-2}
-\sum_{i\in\cE:\atop k\le i\le n-2}\sum_{j\in\cE:\atop k\le j\le i}
+\sum_{i\in\cE:\atop k\le i\le n-k-2}\sum_{j\in\cE:\atop j\le k-2}
+\sum_{i\in\cE:\atop n-k\le i\le n-4}\sum_{j\in\cE:\atop i-k+2\le j\le n-k-2}
\label{computation3}\\
&=
-\sum_{i\in\cE:\atop n-k\le i\le n-2}\sum_{j\in\cE:\atop j\le k-2}
-\sum_{i\in\cE:\atop k\le i\le n-2}\sum_{j\in\cE:\atop k\le j\le i}
+\sum_{i\in\cE:\atop n-k\le i\le n-2}\sum_{j\in\cE:\atop i-k+2\le j\le n-k-2}
\label{computation4}\\
&=
-\sum_{i\in\cE:\atop n-k\le i\le n-2}\sum_{j\in\cE:\atop i-k+2\le j\le i}
-\sum_{i\in\cE:\atop k\le i\le n-2}\sum_{j\in\cE:\atop j\le i-k}
+\sum_{i\in\cE:\atop n-k\le i\le n-2}\sum_{j\in\cE:\atop i-k+2\le j\le n-k-2}
\label{computation5}\\
&=
-\sum_{i\in\cE:\atop k\le i\le n-2}\sum_{j\in\cE:\atop j\le i-k}
-\sum_{i\in\cE:\atop n-k\le i\le n-2}\sum_{j\in\cE:\atop n-k\le j\le i}.
\label{computation6}
\end{align}
We now explain the transformations carried out in the foregoing calculations.
\begin{enumerate}
\item
From~\refp{computation1} to~\refp{computation2}:
the difference of the first and last terms of~\refp{computation1} yields
the first term of~\refp{computation2}. Moreover, we have performed the
substitution $j\mapsto i-j$ in the third term of~\refp{computation1}
($\sum_{j\in\cE:i+k-n+2\le j\le k-2}$) which supplies the third term
of~\refp{computation2} ($\sum_{j\in\cE:i-k+2\le j\le n-k-2}$);
\item
From~\refp{computation2} to~\refp{computation3}:
we have split the sum with respect to $j$ ($j\le i$) of the first term
in~\refp{computation2} into two sums ($\sum_{j\in\cE:j\le k-2}$ and
$\sum_{j\in\cE:k\le j\le i}$).
This gives the two first terms of~\refp{computation3};
\item
From~\refp{computation3} to~\refp{computation4}:
we have simplified the difference of the two sums with respect to $i$
within the first and third terms ($\sum_{i\in\cE:k\le i\le n-2}$ and
$\sum_{i\in\cE:k\le i\le n-k-2}$); this yields the sum with respect to $i$
in the first term in~\refp{computation4} ($\sum_{i\in\cE:n-k\le i\le n-2}$).
Moreover, we have added a null term in the third term: in~\refp{computation4},
$i$ varies from $n-k$ to $n-4$, while in~\refp{computation5}, $i$ varies from
$n-k$ to $n-2$. The sum with respect to $j$ corresponding to the exceeding
index $i=n-2$ actually vanishes;
\item
From~\refp{computation4} to~\refp{computation5}:
we have performed the substitution $j\mapsto i-j$ in the two terms
of~\refp{computation4} ($\sum_{j\in\cE:j\le k-2}$ and $\sum_{j\in\cE:k\le j\le i}$),
which supplies the two first terms of~\refp{computation5}
($\sum_{j\in\cE:i-k+2\le j\le i}$ and $\sum_{j\in\cE:j\le k-i}$);
\item
From~\refp{computation5} to~\refp{computation6}:
the difference of the first and last terms of~\refp{computation5} yields
the last term of~\refp{computation6}.
\end{enumerate}

Finally, plugging~\refp{computation6} into~\refp{intermediate-product}
yields~\refp{lawfree-even-case}.
\fin

\demdeux
As in~\cite{renyi}, formula~\refp{law-product} can be proved by induction.
This way is not constructive but elegant, and we find it interesting to
produce it here.
We formulate the induction hypothesis as follows: for any even $n$,
\\
\textsl{$(\cH_n)$:
For all $m\in\cE$ such that $m\le n$, (for all $j\in\cE$ such that $j\le m$,
$r_{j,m}=r_{j,j}r_{0,m-j}$).}
\\
The initialization $(\cH_0)$ is trivial. Pick now $n\in\cE^*$ and assume that
$\cH_{n-2}$ holds true. Let $k\in\cE$ such that $1\le k\le n-2$.
We have, by $\cH_{n-2}$ and~\refp{rec},
\begin{align*}
r_{k,n}
&=
\frac12\bigg[\sum_{j\in\cE^*:\atop j\le k}\pi_j\,r_{k-j,n-j}
+\sum_{j\in\cE^*:\atop j\le n-k}\pi_j\,r_{k,n-j}\bigg]
\\
&=
\frac12\bigg[\sum_{j\in\cE^*:\atop j\le k}\pi_j\, r_{k-j,k-j}r_{0,n-k}
+\sum_{j\in\cE^*:\atop j\le n-k}\pi_j\,r_{k,k}r_{0,n-j-k}\bigg]
\\
&=
\frac12\bigg[\bigg(\sum_{j\in\cE^*:\atop j\le k}\pi_j\,r_{k-j,k-j}\bigg)r_{0,n-k}
+\bigg(\sum_{j\in\cE^*:\atop j\le n-k}\pi_j\,r_{0,n-j-k}\bigg)r_{k,k}\bigg].
\end{align*}
Remarking, by the Markov property, that
$$
\frac12\,\pi_j\,r_{k-j,k-j}=\P\{S_1>0,\tau_0=j,T_k=k\}\quad\mbox{and}\quad
\frac12\,\pi_j\,r_{0,n-j-k}=\P\{S_1<0,\tau_0=j,T_k=k\},
$$
we obtain
\beq
r_{k,n}=\P\{\tau_0\le k,T_k=k\}r_{0,n-k}+\P\{\tau_0\le n-k,T_{n-k}=0\}r_{k,k}.
\label{intermediate-product2}
\eeq
Now, we have for even $k$
\begin{align}
\P\{\tau_0\le k,T_k=k\}
&=
\P\{T_k=k\}-\P\{\tau_0>k,T_k=k\}
\nonumber\\
&=
\P\{S_1\ge 0,\dots,S_k\ge 0\}-\P\{S_1>0,\dots,S_k>0\}
\nonumber\\
&=
\P\{\tau_{-1}>k\}-p\,\P_1\{\tau_0>k-1\}
=\P\{\tau_{-1}>k\}-p\,\P\{\tau_{-1}>k-1\}
\nonumber\\
&=
\P\{\tau_{-1}>k\}-p\,\P\{\tau_{-1}>k\}
=q\,\P\{\tau_{-1}>k\}=q\,\P\{T_k=k\}
\nonumber\\
&=
q\,r_{k,k}.
\label{first-termP}
\end{align}
Similarly, we have
\begin{align}
\P\{\tau_0\le n-k,T_{n-k}=0\}
&=
p\,\P\{\tau_1>n-k\}=p\,\P\{T_{n-k}=0\}
\nonumber\\
&=
p\,r_{0,n-k}.
\label{second-termP}
\end{align}
Finally, putting~\refp{first-termP} and~\refp{second-termP}
into~\refp{intermediate-product2}, we derive~\refp{law-product}.
\fin

\demtrois
We propose a last and very short proof of~\refp{law-product}.
The probabilities $r_{k,n}$, $k,n\in\cE$, are generated by the function $\tG$
which was introduced in Subsection~\ref{section-free}.
By~\refp{geneGfree-even-product}, the quantity $\tG(x,y)$ can be factorized into
the product of $\tG(x,0)$ and $\tG(0,y)$. Writing that
$$
\tG(x,0)\tG(0,y)=\sum_{k,l\in\cE}(r_{k,k} r_{0,l})x^ky^l
=\sum_{k,n\in\cE:\atop k\le n}(r_{k,k} r_{0,n-k})x^ky^{n-k},
$$
we conclude, by identification, that~\refp{law-product} holds true.
\fin

%
%%%%%%%%%%%%%%
\bco{
In the symmetric case ($p=q=1/2$), the well-known following expression holds
for even integers $n,k$ such that $0\le k\le n$:
$$
r_{k,n}=\frac{1}{2^n}\binom{k}{k/2}\binom{(n-k)}{(n-k)/2}.
$$
}
%%%%%%%%%%%%%%
%
\dem
In the case where $p=1/2$, we have the particular identities
$$
\P\{\tau_0=j\}=\P\{S_{j-2}=0\}-\P\{S_j=0\} \quad\mbox{and}\quad
\P\{\tau_{-1}=j\}=\P\{\tau_0=j+1\}
$$
which are due to the fact that $b_{j-2}-b_j/4=a_{j-2}$ and to~\refp{law-tau1}
respectively, and then
\begin{align*}
r_{k,k}
&=
\P\{\tau_{-1}\ge k+1\}=\sum_{j\in\cO:\atop j\ge k+1}\P\{\tau_{-1}=j\}
=\sum_{j\in\cO:\atop j\ge k+1}\P\{\tau_0=j+1\}
=\sum_{j\in\cE:\atop j\ge k+2}\P\{\tau_0=j\}
\\
&=
\sum_{j\in\cE:\atop j\ge k+2}[\P\{S_{j-2}=0\}-\P\{S_j=0\}]
=\P\{S_k=0\}=\frac{1}{2^k}\binom{k}{k/2}.
\end{align*}
Analogously,
$$
r_{0,n-k}=\P\{\tau_1\ge n-k+1\}=\P\{S_{n-k}=0\}=\frac{1}{2^{n-k}}
\binom{(n-k)}{(n-k)/2}.
$$
We conclude with the help of~\refp{law-product}.
\fin

%
%%%%%%%%%%%%%%
\brem{
Let us compute the limit of $r_{k,n}$ as $n$ tends to $\infty$.
We have to distinguish the cases $n\in\cE$ and $n\in\cO$.
On one hand, by~\refp{limit-r0n} and~\refp{law-product}, we have
$$
\lim_{n\to \infty\atop  n\in\cE} r_{k,n}
=\ind_{\cE}(k) \lim_{n\to \infty\atop  n\in\cE} r_{k,k}r_{0,n-k}
=\ind_{\cE}(k) r_{k,k}\lim_{n\to \infty\atop  n\in\cE} r_{0,n}
=\ind_{\cE}(k) r_{k,k}\bigg(2-\frac 1q\bigg)^+.
$$
On the other hand, since
\begin{align*}
\lim_{n\to \infty} \bigg[2p\sum_{i\in\cE:\atop n-k\le i\le n} a_i (pq)^{i/2}
-4\sum_{i\in\cE:\atop n-k\le i\le n-2}
\bigg(\sum_{j\in \cE:\atop j\le i+k-n} a_ja_{i-j}\bigg) (pq)^{i/2+1}\bigg]
=0
\end{align*}
we can see from~\refp{lawfree-odd-case} that for odd $k$
$$
\lim_{n\to \infty\atop  n\in\cO} r_{k,n}=0.
$$
Assume now that $k$ is even. By~\refp{lawfree}, we get
\begin{align*}
\lim_{n\to \infty\atop  n\in\cO} r_{k,n}
&=
\lim_{n\to \infty\atop  n\in\cO}
\bigg[2p\sum_{i\in\cE:\atop n-k\le i\le n} a_i (pq)^{i/2}
-4\sum_{i\in\cE:\atop n-k\le i\le n-2}
\bigg(\sum_{j\in \cE:\atop j\le i+k-n} a_ja_{i-j}\bigg) (pq)^{i/2+1}
\\
&\hphantom{=\,}
+2q\sum_{i\in\cE:\atop k\le i\le n} a_i (pq)^{i/2}
-4\sum_{i\in\cE:\atop k\le i\le n-2}
\bigg(\sum_{j\in \cE:\atop j\le i-k} a_ja_{i-j}\bigg) (pq)^{i/2+1}\bigg]
\\
&=
\lim_{n\to \infty\atop  n\in\cO}
\bigg[2p\sum_{i\in\cE:\atop n+1-k\le i\le n-1} a_i (pq)^{i/2}
-4\sum_{i\in\cE:\atop n+1-k\le i\le n-3}
\bigg(\sum_{j\in \cE:\atop j\le i+k-n-1} a_ja_{i-j}\bigg) (pq)^{i/2+1}
\\
&\hphantom{=\,}
+2q\sum_{i\in\cE:\atop k\le i\le n-1} a_i (pq)^{i/2}
-4\sum_{i\in\cE:\atop k\le i\le n-3}
\bigg(\sum_{j\in \cE:\atop j\le i-k} a_ja_{i-j}\bigg) (pq)^{i/2+1}\bigg]
\\
\lqn{}
&=
\lim_{n\to \infty\atop  n\in\cO}
\bigg[2p\sum_{i\in\cE:\atop n+1-k\le i\le n+1} a_i (pq)^{i/2}
-4\sum_{i\in\cE:\atop n+1-k\le i\le n-1}
\bigg(\sum_{j\in \cE:\atop j\le i+k-n-1} a_ja_{i-j}\bigg) (pq)^{i/2+1}
\\
&\hphantom{=\,}
+2q\sum_{i\in\cE:\atop k\le i\le n+1} a_i (pq)^{i/2}
-4\sum_{i\in\cE:\atop k\le i\le n-1}
\bigg(\sum_{j\in \cE:\atop j\le i-k} a_ja_{i-j}\bigg) (pq)^{i/2+1}
\\
&\hphantom{=\,}
-2p\,a_{n+1} (pq)^{(n+1)/2}
+4 \bigg(\sum_{j\in \cE:\atop j\le k-2} a_ja_{n-1-j}\bigg) (pq)^{(n+1)/2}
\\
&\hphantom{=\,}
-2q\,a_{n+1} (pq)^{(n+1)/2}
+4 \bigg(\sum_{j\in \cE:\atop j\le n-1-k} a_ja_{n-1-j}\bigg) (pq)^{(n+1)/2}\bigg].
\end{align*}
The last four terms within the above limit tend to 0 and it remains the first
four terms which are nothing but $r_{k,n+1}$. Hence,
$$
\lim_{n\to \infty\atop  n\in\cO} r_{k,n}
=\lim_{n\to \infty\atop  n\in\cO} r_{k,n+1}
=\lim_{n\to \infty\atop  n\in\cE} r_{k,n}.
$$
As a result, the foregoing discussion shows that
$$
\lim_{n\to \infty\atop  n\in\cE} r_{k,n}
=\lim_{n\to \infty\atop  n\in\cO} r_{k,n}
=\lim_{n\to \infty} r_{k,n}
$$
which supplies the distribution of the total sojourn time of the random walk
in $\Z^+$:
$$
\P\{T_{\infty}=k\}=\ind_{\cE}(k) r_{k,k}\bigg(2-\frac 1q\bigg)^+.
$$
By~\refp{lawfree-part-cases} and~\refp{elembis}, we can rewrite $r_{k,k}$ as
\begin{align*}
r_{k,k}
&=
1-2q\bigg[\sum_{i\in\cE} a_i (pq)^{i/2}
-\sum_{i\in\cE:\atop i\ge k} a_i (pq)^{i/2}\bigg]
=1-2q\bigg[\frac{1-|p-q|}{4pq}-\sum_{i\in\cE:\atop i\ge k} a_i (pq)^{i/2}\bigg]
\\
&=
\bigg(2-\frac 1p\bigg)^+ +2q\sum_{i\in\cE:\atop i\ge k} a_i (pq)^{i/2}.
\end{align*}
With this representation of $r_{k,k}$ at hand, we can compute the
probability $\P\{T_{\infty}<\infty\}$. First, if $p\ge 1/2$, we plainly
have $\P\{T_{\infty}=k\}=0$ for all $k$ and then $\P\{T_{\infty}<\infty\}=0$.
Second, if $p<1/2$,
\begin{align*}
\P\{T_{\infty}<\infty\}
&=
\bigg(2-\frac 1q\bigg)\sum_{k\in\cE} r_{k,k}
=2(2q-1)\sum_{k\in\cE}\bigg(\sum_{i\in\cE:\atop i\ge k} a_i (pq)^{i/2}\bigg)
\\
&=
2(2q-1)\sum_{i\in\cE} \bigg(\sum_{k\in\cE:\atop k\le i} 1\bigg)a_i (pq)^{i/2}
=(2q-1)\sum_{i\in\cE} (i+2)a_i (pq)^{i/2}
\\
&=
(2q-1)\sum_{i\in\cE} b_i (pq)^{i/2}
=\frac{q-p}{A(pq)}=1.
\end{align*}
The next to last equality comes from~\refp{elem}. Therefore we have obtained the result
$$
\P\{T_{\infty}<\infty\}=\left\{\begin{array}{ll}
0 & \mbox{if $p\ge 1/2$,}
\\
1 & \mbox{if $p<1/2$.}
\end{array}\right.
$$
which is in good agreement with the transient or recurrent character of the random walk.
Indeed, if $p>1/2$, $\lim_{n\to \infty}S_n=+\infty$ a.s.;
if $p<1/2$, $\lim_{n\to \infty}S_n=-\infty$ a.s.;
if $p=1/2$, $\limsup_{n\to \infty}S_n=+\infty$ a.s.
Hence, one has $T_{\infty}=\infty$ a.s. for $p\ge 1/2$, while
$T_{\infty}<\infty$ a.s. for $p<1/2$.
}
%%%%%%%%%%%%%%
%
%%%%%%%%%%%%%%
\bpr{
For odd integers $n,k$ such that $0\le k\le n-1$, the following relationship
holds:
\beq
r_{k,n}+r_{k+1,n}=r_{k+1,n+1}.
\label{rec2}
\eeq
}
%%%%%%%%%%%%%%
%
\demun
We first derive~\refp{rec2} by using the explicit results obtained in
Theorem~\ref{theorem1bis}.
On one hand, since $k$ and $n$ are odd, \refp{lawfree-odd-case} yields
\begin{align*}
r_{k,n}+r_{k+1,n}
&=
2p\sum_{i\in\cE:\atop n-k\le i\le n-1} a_i (pq)^{i/2}
-4\sum_{i\in\cE:\atop n-k\le i\le n-3}
\bigg(\sum_{j\in \cE:\atop j\le i+k-n} a_ja_{i-j}\bigg) (pq)^{i/2+1}
\\
&\hphantom{=\,}
+2q\sum_{i\in\cE:\atop k+1\le i\le n-1} a_i (pq)^{i/2}
-4\sum_{i\in\cE:\atop k+1\le i\le n-3}
\bigg(\sum_{j\in \cE:\atop j\le i-k-1} a_ja_{i-j}\bigg) (pq)^{i/2+1}.
\end{align*}
On the other hand, since $k+1$ and $n+1$ are even, \refp{lawfree-even-case}
yields
\begin{align*}
r_{k+1,n+1}
&=
2p\sum_{i\in\cE:\atop n-k\le i\le n+1} a_i (pq)^{i/2}
-4\sum_{i\in\cE:\atop n-k\le i\le n-1}
\bigg(\sum_{j\in \cE:\atop j\le i+k-n} a_ja_{i-j}\bigg) (pq)^{i/2+1}
\\
&\hphantom{=\,}
+2q\sum_{i\in\cE:\atop k+1\le i\le n+1} a_i (pq)^{i/2}
-4\sum_{i\in\cE:\atop k+1\le i\le n-1}
\bigg(\sum_{j\in \cE:\atop j\le i-k-1} a_ja_{i-j}\bigg) (pq)^{i/2+1}.
\end{align*}
We then deduce, by remarking that
$\sum_{j\in \cE:\atop j\le k-1} a_ja_{n-1-j}
=\sum_{j\in \cE:\atop n-k\le j\le n-1} a_ja_{n-1-j}$
and by using~\refp{convol}, that
\begin{align*}
\lqn{r_{k+1,n+1}-(r_{k,n}+r_{k+1,n})}\\[-3ex]
&=
(2p+2q)a_{n+1} (pq)^{(n+1)/2}
\\
&\hphantom{=\,}
-4\bigg(\sum_{j\in \cE:\atop j\le k-1} a_ja_{n-1-j}\bigg)(pq)^{(n-1)/2+1}
-4\bigg(\sum_{j\in \cE:\atop j\le n-k-2} a_ja_{n-1-j}\bigg)(pq)^{(n-1)/2+1}
\\
\lqn{}
&=
2\bigg(a_{n+1}-2\sum_{j\in \cE:\atop j\le n-1} a_ja_{n-1-j}\bigg)(pq)^{(n+1)/2}=0.
\end{align*}
The relationship~\refp{rec2} is checked.
\fin

\demdeux
We can prove~\refp{rec2} directly without using any explicit expression.
Pick two odd integers $n,k$ such that $0\le k\le n-1$. We have
\beq
r_{k,n}+r_{k+1,n}=\P\{T_n=k\}+\P\{T_n=k+1\}=\P\{T_n\in\{k,k+1\}\}.
\label{intermediate-rec2}
\eeq
Let us introduce the last passage time by 0, say $\s_0$, for the random walk.
Is is plain that $\s_0$ is even and that there is an even number of $\d_j$ up to
time $\s_0$ which are equal to one. Afterwards, it remains an odd number
(this is $n-\s_0$) of steps up to time $n$.
One has either $S_{\s_0+1}>0,S_{\s_0+2}>0,\dots,S_n>0$ or
$S_{\s_0+1}<0,S_{\s_0+2}<0,\dots,S_n<0$.
The corresponding $\d_j$, $\s_0+1\le j\le n$,
are either all equal to one or all equal to zero. More precisely,
\bitem
\item
if $T_n=k$, since $k$ is odd, one has in this case
$S_{\s_0+1}>0,S_{\s_0+2}>0,\dots,S_n>0$.
Then, necessarily, $S_{n+1}\ge 0$ which entails that $\d_{n+1}=1$ and
$T_{n+1}=T_n+1=k+1$;
\item
if $T_n=k+1$, since $k+1$ is even, one has in this case
$S_{\s_0+1}<0,S_{\s_0+2}<0,\dots,S_n<0$.
Then, necessarily, $S_{n+1}\le 0$ which entails that $\d_{n+1}=0$
and $T_{n+1}=T_n=k+1$.
\eitem
This discussion implies the inclusion
$\{T_n\in\{k,k+1\}\}\subset\{T_{n+1}=k+1\}.$
Conversely, since $T_{n+1}-T_n=\d_n\in\{0,1\}$, the equality $T_{n+1}=k+1$
implies $T_n\in\{k,k+1\}$ which proves the inclusion
$\{T_{n+1}=k+1\}\subset\{T_n\in\{k,k+1\}\}.$ As a byproduct,
the equality $\{T_n\in\{k,k+1\}\}=\{T_{n+1}=k+1\}$ holds and this proves,
referring to~\refp{intermediate-rec2}, the relationship~\refp{rec2}.
\fin

%%%%%%%%%%%%%%%%%%%%%%%%%%%%%%%%%%%%%%%%%%%%%%%%%%%%%%%%%%%%%%%%%%%%%%%%%%%%%%%
%%%%%%%%%%%%%%%%%%%%%%%%%        Bridge        %%%%%%%%%%%%%%%%%%%%%%%%%%%%%%%%
%%%%%%%%%%%%%%%%%%%%%%%%%%%%%%%%%%%%%%%%%%%%%%%%%%%%%%%%%%%%%%%%%%%%%%%%%%%%%%%
\section{The case $F=\{0\}$ (bridge of the random walk)}\label{section-bridge}

In this part, we shall set for simplifying $r_{k,n}^{\{0\}}=r_{k,n}^{0}$,
$H^{\{0\}}(z)=H^{0}(z)$ and $G^{\{0\}}(x,y)=G^{0}(x,y)$.
We consider the distribution of the sojourn time $T_n$ subject to the
condition that $S_n=0$, that is we are dealing with the so-called bridge of
the random walk pinned at zero at times $0$ and $n$.
The condition $S_n=0$ can be fulfilled only when $n$ is even. So, we make
here the assumption that $n$ is an even integer throughout this section.

%%%%%%%%%%%%%%%%%%%%%%%%%%%%%%%%%%%%%%%%%%%%%%%%%%%%%%%%%%%%%%%%%%%%%%%%%%%%%%%
%%%%%%%%%%%%%%%%%%%%%%%%%  Generating function   %%%%%%%%%%%%%%%%%%%%%%%%%%%%%%
%%%%%%%%%%%%%%%%%%%%%%%%%%%%%%%%%%%%%%%%%%%%%%%%%%%%%%%%%%%%%%%%%%%%%%%%%%%%%%%
\subsection{Generating function}

In the light of Theorem~\ref{th-generalG}, we need to evaluate the functions
$H_i^0(z)$, $G^0(x,0)$ and $G^0(0,y)$.
We have $H_i^0(z)=H_{-i}(z)$ and then, by~\refp{geneH} and~\refp{GF(x,0)},
\begin{align*}
G^0(x,0)
&=
H_0^0(x)-\frac{1-\bA(x)}{2px}\,H_{-1}^0(x)
=H_0(x)-\frac{1-\bA(x)}{2px}\,H_1(x)
\\
&=
\frac{1}{\bA(x)}\bigg[1-\frac{1-\bA(x)}{2px}\,\frac{2px}{1+\bA(x)}\bigg]
=\frac{2}{1+\bA(x)}.
\end{align*}
Similarly,
$$
G^0(0,y)=\frac{2}{1+\bA(y)}.
$$
From this and~\refp{generalG}, we immediately extract $G^0(x,y)$.
%
%%%%%%%%%%%%%%
\bth{\label{theorem2}
The generating function $G^0$ is given by
\beq
G^0(x,y)=\frac{2}{\bA(x)+\bA(y)}
\label{geneGbridge}
\eeq
where $\bA(z)=\sqrt{1-4pqz^2}$.
}
%%%%%%%%%%%%%%
%
The result~\refp{geneGbridge} is remarkably simple.

%%%%%%%%%%%%%%%%%%%%%%%%%%%%%%%%%%%%%%%%%%%%%%%%%%%%%%%%%%%%%%%%%%%%%%%%%%%%%%%
%%%%%%%%%%%%%%%%%%%%%%%%%         Inverting      %%%%%%%%%%%%%%%%%%%%%%%%%%%%%%
%%%%%%%%%%%%%%%%%%%%%%%%%%%%%%%%%%%%%%%%%%%%%%%%%%%%%%%%%%%%%%%%%%%%%%%%%%%%%%%
\subsection{Distribution of the sojourn time}

%
%%%%%%%%%%%%%%
\bth{\label{theorem2bis}
Assume that $n$ is even.
The probability $r_{k,n}^0=\P\{T_n=k,S_n=0\}$ admits the following expression:
\beq
r_{k,n}^0=\left\{\begin{array}{ll}
\dis\frac{2}{n+2}\binom{n}{n/2} (pq)^{n/2}
& \mbox{if $k$ is even such that $0\le k\le n$,}
\\[2ex]
0 & \mbox{if $k$ is odd or such that $k\ge n+1$.}
\end{array}\right.
\label{law-bridge}
\eeq
The distribution of $T_n$ for the bridge of the random walk is given by
$$
\P\{T_n=k\,|\,S_n=0\}=\left\{\begin{array}{ll}
\dis\frac{2}{n+2} & \mbox{if $k$ is even such that $0\le k\le n$,}
\\[2ex]
0 & \mbox{if $k$ is odd or such that $k\ge n+1$.}
\end{array}\right.
$$
that is, the conditioned random variable $(T_n=k\,|\,S_n=0)$ is
uniformly distributed on the set $\{0,2,4,\dots,n-2,n\}$.
}
%%%%%%%%%%%%%%
%
\demun
By~\refp{1surA+A}, we rewrite $G^0(x,y)$ as
$$
G^0(x,y)=2\sum_{i,j\in\cE}a_{i+j} (pq)^{(i+j)/2}x^iy^j
=2\sum_{k,n\in\cE:\atop k\le n}a_n (pq)^{n/2}x^ky^{n-k}.
$$
From this, we immediately extract the announced expression for $r_{k,n}^0$.
Moreover, we plainly have
$$
\P\{T_n=k\,|\,S_n=0\}=\frac{\P\{T_n=k,S_n=0\}}{\P\{S_n=0\}}
$$
with
$$
\P\{S_n=0\}=\binom{n}{n/2} (pq)^{n/2}
$$
and this prove the assertion related to the uniform law.
\fin

\demdeux
We find it interesting to give a proof by induction relying on the recurrent
relationship~\refp{rec} which actually holds even in the present case, namely:
\beq
r_{k,n}^0=\frac12\bigg[\sum_{j\in\cE^*:\atop j\le k}\pi_j\,r_{k-j,n-j}^0
+\sum_{j\in\cE^*:\atop j\le n-k}\pi_j\,r_{k,n-j}^0\bigg].
\label{rec-bridge}
\eeq
Indeed, the proof of~\refp{rec} may be adapted \textit{mutatis mutandis}
in order to derive~\refp{rec-bridge}.
We formulate the induction hypothesis as follows:
\\
\textsl{$(\cH_n^0)$:
For all $m\in\cE$ such that $m\le n$, (for all $j\in\cE$ such that $j\le m$,
$r_{j,m}^0=r_{0,m}^0=2a_m(pq)^{m/2}$).}
\\
The initialization $(\cH_0^0)$ is trivial. Pick $n\in\cE^*$ and assume that
$\cH_{n-2}^0$ holds true. Let $k\in\cE$ such that $1\le k\le n-1$.
We have, by $\cH_{n-2}^0$, \refp{convol}, \refp{law-tau0} and~\refp{rec-bridge},
\begin{align*}
r_{k,n}^0
&=
\frac12\bigg[\sum_{j\in\cE^*:\atop j\le k}\pi_j\,r_{0,n-j}^0
+\sum_{j\in\cE^*:\atop j\le n-k}\pi_j\,r_{0,n-j}^0\bigg]
=4\bigg[\sum_{j\in\cE^*:\atop j\le k}a_{j-2} a_{n-j}
+\sum_{j\in\cE^*:\atop j\le n-k}a_{j-2} a_{n-j}\bigg] (pq)^{n/2}
\\
&=
4\bigg[\sum_{j\in\cE:\atop j\le k-2}a_j a_{n-2-j}
+\sum_{j\in\cE:\atop k\le j\le n-2}a_j a_{n-2-j}\bigg] (pq)^{n/2}
=4\bigg[\sum_{j\in\cE:\atop j\le n-2}a_j a_{n-2-j}\bigg] (pq)^{n/2}
\\
&=
2a_n(pq)^{n/2}=r_{0,n}^0.
\end{align*}
It remains to see that this formula also holds for $k=n$. By~\refp{convol},
\refp{law-Sk-part}, \refp{law-tau0} and~\refp{part-cases}, we have
\begin{align*}
r_{n,n}^0
&=
\P\{S_n=0\}-\frac{1}{2p}\sum_{j\in\cO:\atop j\le n-1} \pi_{j+1}\P_{-1}\{S_{n-j}=0\}
\\
&=
(n+2)a_n(pq)^{n/2}-\sum_{j\in\cE:\atop j\le n-2} a_jb_{n-j}(pq)^{n/2}
\\
&=
[(n+3)a_n-\frac14 \,b_{n+2}](pq)^{n/2}
=[(n+3)a_n-(n+1)a_n](pq)^{n/2}
\\
&=
2a_n(pq)^{n/2}.
\end{align*}
By~\refp{part-cases}, we also have
$$
r_{0,n}^0=\P\{S_n=0\}-\frac{1}{2q}\sum_{j\in\cO:\atop j\le n-1}
\pi_{j+1}\P_1\{S_{n-j}=0\}.
$$
Because of the ``duality'' relationship
$q\,\P_1\{S_{n-j}=0\}=p\,\P_{-1}\{S_{n-j}=0\}$,
we immediately see that $r_{n,n}^0=r_{0,n}^0.$
The validity of $(\cH_n^0)$ is acquired
and the proof of~\refp{law-bridge} is finished.
\fin

%%%%%%%%%%%%%%%%%%%%%%%%%%%%%%%%%%%%%%%%%%%%%%%%%%%%%%%%%%%%%%%%%%%%%%%%%%%%%%%
%%%%%%%%%%%%%%%%%%%%%%%%%        Meander       %%%%%%%%%%%%%%%%%%%%%%%%%%%%%%%%
%%%%%%%%%%%%%%%%%%%%%%%%%%%%%%%%%%%%%%%%%%%%%%%%%%%%%%%%%%%%%%%%%%%%%%%%%%%%%%%
\section{The cases $F=\Z^{+*}$ and $F=\Z^{-*}$}

In this part, we shall set for simplifying $r_{k,n}^{\Z^{+*}}=r_{k,n}^+$,
$G^{\Z^{+*}}(x,y)=G^+(x,y)$ and $H^{\Z^{+*}}(z)=H^+(z)$. We shall also
use similar notations with minus signs for the study of the case $F=\Z^{-*}$.
Some expressions for the $r_{k,n}^+$'s can be obtained from
Theorem~\ref{theorem0bis} by summing the $r_{k,n}^j$, $j\in\Z^{+*}$ or
$j\in\Z^{-*}$. Since $r_{k,n}^j=0$ for $|j|>n$, the corresponding sums reduce to
$$
r_{k,n}^+=\sum_{j=1}^n r_{k,n}^j\quad\mbox{and}\quad
r_{k,n}^-=\sum_{j=-n}^{-1} r_{k,n}^j.
$$
We propose other expressions which can be deduced from the generating functions
$G^+$ and $G^-$.

%%%%%%%%%%%%%%%%%%%%%%%%%%%%%%%%%%%%%%%%%%%%%%%%%%%%%%%%%%%%%%%%%%%%%%%%%%%%%%%
%%%%%%%%%%%%%%%%%%%%%%%%%  Generating function   %%%%%%%%%%%%%%%%%%%%%%%%%%%%%%
%%%%%%%%%%%%%%%%%%%%%%%%%%%%%%%%%%%%%%%%%%%%%%%%%%%%%%%%%%%%%%%%%%%%%%%%%%%%%%%
\subsection{Generating function}

We first consider the case where $F=\Z^{+*}$.
We need to evaluate the functions $H_i^+(z)$ for $i\in\{-1,0,1\}$, $G^+(x,0)$
and $G^+(0,y)$. On one hand, we have for $i\in\{-1,0,1\}$ (and then $i\ge 0$)

\begin{align*}
H_i^+(z)
&=
\sum_{j\in\Z^{+*}-i}H_j(z)
=\frac{1}{\bA(z)}\sum_{j=1-i}^{\infty} \bigg(\frac{1-\bA(z)}{2qz}\bigg)^j
\\
&=
\frac{1}{\bA(z)(1-\frac{1-\bA(z)}{2qz})}\bigg(\frac{1-\bA(z)}{2qz}\bigg)^{1-i}
\\
&=
\frac{1-\bA(z)}{\bA(z)[2qz-1+\bA(z)]}\bigg(\frac{1-\bA(z)}{2qz}\bigg)^{-i}.
\end{align*}
Invoking~\refp{geneH} and~\refp{identityA}, we find that
$$
H_i^+(z)=\frac{2pz-1+\bA(z)}{2(1-z)\bA(z)}\bigg(\frac{1-\bA(z)}{2qz}\bigg)^{-i}
=\frac{2pz-1+\bA(z)}{2(1-z)\bA(z)}\bigg(\frac{1+\bA(z)}{2pz}\bigg)^i.
$$
On the other hand, by~\refp{GF(x,0)},
\begin{align*}
G^+(x,0)
&=
H_0^+(x)-\frac{1-\bA(x)}{2px}\,H_{-1}^+(x)
=H_0^+(x)\bigg[1-\frac{1-\bA(x)}{2px}\,\frac{2px}{1+\bA(x)}\bigg]
\\
&=
\frac{2\bA(x)}{1+\bA(x)}\,H_0^+(x)
=\frac{2px-1+\bA(x)}{(1-x)[1+\bA(x)]}.
\end{align*}
Similarly,
$$
G^+(0,y)=H_0^+(y)-\frac{1-\bA(y)}{2qy}\,H_1^+(y)
=H_0^+(y)\bigg[1-\frac{1-\bA(y)}{2qy}\,\frac{2qy}{1-\bA(y)}\bigg]=0.
$$
With this at hand, we can derive $G^+(x,y)$. Indeed, using the general
formula~\refp{generalG}, we obtain
$$
G^+(x,y)=\frac{1+\bA(x)}{\bA(x)+\bA(y)}\,G^+(x,0)
=\frac{2px-1+\bA(x)}{(1-x)[\bA(x)+\bA(y)]}.
$$
Exactly in the same way, we could find the result related to the case
$F=\Z^{-*}$. We state both results in the theorem below.
%
%%%%%%%%%%%%%%
\bth{\label{theorem3}
The generating functions $G^+$ and $G^-$ are given by
\beq
G^+(x,y)=\frac{2px-1+\bA(x)}{(1-x)[\bA(x)+\bA(y)]}
\quad\mbox{and}\quad
G^-(x,y)=\frac{2qy-1+\bA(y)}{(1-y)[\bA(x)+\bA(y)]}
\label{geneGpositive}
\eeq
where $\bA(z)=\sqrt{1-4pqz^2}$.
}
%%%%%%%%%%%%%%
%
%%%%%%%%%%%%%%
\brem{\label{rem-sum}
We have the following relationship between the generating functions
$G$, $G^+$, $G^-$, $G^0$:
$$
G^+(x,y)+G^-(x,y)+G^0(x,y)=G(x,y).
$$
Indeed, it is clearly due to the fact that
$\{S_n\in\Z^{+*}\}\cup\{S_n\in\Z^{-*}\}\cup\{S_n=0\}=\Omega$
which implies that
$$
r_{k,n}^++r_{k,n}^-+r_{k,n}^0=r_{k,n}.
$$
As a checking, we propose to compute directly the sum $G^++G^-+G^0$
by using the obtained expressions~\refp{geneGfree}, \refp{geneGbridge}
and~\refp{geneGpositive}:
$$
G^+(x,y)+G^-(x,y)+G^0(x,y)
=\frac{1}{\bA(x)+\bA(y)}\bigg[\frac{2px-1+\bA(x)}{1-x}
+\frac{2qy-1+\bA(y)}{1-y}+2\bigg].
$$
The term within brackets can be written as
$$
\frac{2px-1}{1-x}+\frac{2qy-1}{1-y}+2
=\bigg[\frac{2px-1}{1-x}+2p\bigg]+\bigg[\frac{2qy-1}{1-y}+2q\bigg]
=\frac{2p-1}{1-x}+\frac{2q-1}{1-y}
$$
and then
$$
G^+(x,y)+G^-(x,y)+G^0(x,y)
=\frac{1}{\bA(x)+\bA(y)}\bigg[\frac{2p-1+\bA(x)}{1-x}
+\frac{2q-1+\bA(y)}{1-y}\bigg]
=G(x,y).
$$
}

%%%%%%%%%%%%%%%%%%%%%%%%%%%%%%%%%%%%%%%%%%%%%%%%%%%%%%%%%%%%%%%%%%%%%%%%%%%%%%%
%%%%%%%%%%%%%%%%%%%%%%%%%         Inverting      %%%%%%%%%%%%%%%%%%%%%%%%%%%%%%
%%%%%%%%%%%%%%%%%%%%%%%%%%%%%%%%%%%%%%%%%%%%%%%%%%%%%%%%%%%%%%%%%%%%%%%%%%%%%%%
\subsection{Distribution of the sojourn time}

From Theorem~\ref{theorem3}, we derive the coefficients $r_{k,n}^+$ and $r_{k,n}^-$.
%
%%%%%%%%%%%%%%
\bth{\label{theorem3bis}
The probabilities $r_{k,n}^+=\P\{T_n=k,S_n>0\}$ and $r_{k,n}^-=\P\{T_n=k,S_n<0\}$
admit the following expressions: for $0\le k\le n$,
\begin{align}
r_{k,n}^+
&=
\ind_{\cE}(n-k)\bigg[2p\sum_{i\in\cE:\atop n-k\le i\le n-1} a_i (pq)^{i/2}
-4\sum_{i\in\cE:\atop n-k\le i\le n-2}
\bigg(\sum_{j\in \cE:\atop j\le i+k-n} a_ja_{i-j}\bigg) (pq)^{i/2+1}\bigg],
\label{law-positive}
\\
r_{k,n}^-
&=
\ind_{\cE}(k)\bigg[2q\sum_{i\in\cE:\atop k\le i\le n-1} a_i (pq)^{i/2}
-4\sum_{i\in\cE:\atop k\le i\le n-2}
\bigg(\sum_{j\in \cE:\atop j\le i-k} a_ja_{i-j}\bigg) (pq)^{i/2+1}\bigg],
\label{law-negative}
\end{align}
where $a_i=\dis\frac{1}{i+2}\binom{i}{i/2}$ for $i\in\cE$.
}
%%%%%%%%%%%%%%
%
\dem
We invert the generating function $G^+$. For its expansion, we refer to the
proof of Theorem~\ref{theorem1bis}. We have
\begin{align*}
G^+(x,y)
&=
\frac{2px}{1-x}\,\frac{1}{\bA(x)+\bA(y)}
-\frac{1-\bA(x)}{1-x}\,\frac{1}{\bA(x)+\bA(y)}
\\
&=
2p\sum_{k,n\in\N:\atop k\le n,} \ind_{\cE}(n-k)
\bigg[\sum_{i\in\cE:\atop n-k\le i\le n} a_i (pq)^{i/2}\bigg] x^{k+1}y^{n-k}
\\
&\hphantom{=\,}
-4\sum_{k,n\in\N:\atop k\le n} \ind_{\cE}(n-k)\bigg[\sum_{i\in\cE:\atop i\le n-2}
\bigg(\sum_{j\in \cE:\atop j\le i+k-n} a_ja_{i-j}\bigg) (pq)^{i/2+1}\bigg] x^ky^{n-k}.
\end{align*}
Performing the substitution $(k,n)\mapsto (k-1,n-1)$ in the first term of the
last above equality, we get
\begin{align*}
G^+(x,y)
&=
2p\sum_{k,n\in\N:\atop k\le n,} \ind_{\cE}(n-k)
\bigg[\sum_{i\in\cE:\atop n-k\le i\le n-1} a_i (pq)^{i/2}\bigg] x^ky^{n-k}
\\
&\hphantom{=\,}
-4\sum_{k,n\in\N:\atop k\le n} \ind_{\cE}(n-k)\bigg[\sum_{i\in\cE:\atop i\le n-2}
\bigg(\sum_{j\in \cE:\atop j\le i+k-n} a_ja_{i-j}\bigg) (pq)^{i/2+1}\bigg] x^ky^{n-k}.
\end{align*}
Formula~\refp{law-positive} ensues by identification.
Formula~\refp{law-negative} can be deduced from~\refp{law-positive} by
invoking the duality argument mentioned in Remark~\ref{rem-duality}: it
suffices to interchange $p$ and $q$ on one hand, and $k$ and $n-k$ on the
other hand.
\fin

%
%%%%%%%%%%%%%%
\brem{
By comparing~\refp{lawfree} and~\refp{law-positive}, we can see that, for odd
integer $n$, $r_{k,n}^+=r_{k,n}$ if $k$ is odd, $r_{k,n}^+=0$ if $k$ is even, and,
for even integer $n$, $r_{0,n}^+=r_{1,n}^+=0$ and $r_{2,n}^+=p^2r_{0,n-2}^0$.
These relations can be directly checked. For instance, in the last case,
we can easily observe that the conditions $T_n=2$ and $S_n>0$ are fulfilled
only in the case where $S_0\le 0,S_1\le 0,\dots,S_{n-3}\le 0,$
$S_{n-2}=0,S_{n-1}=1$ and $S_n=2$. Thus,
\begin{align*}
\P\{T_n=2,S_n>0\}
&=
\P\{S_0\le 0,S_1\le 0,\dots,S_{n-3}\le 0,S_{n-2}=0,X_{n-1}=X_n=1\}
\\
&=
p^2\P\{T_{n-2}=0,S_{n-2}=0\}
\end{align*}
which is nothing but $r_{2,n}^+=p^2r_{0,n-2}^0$.
}
%%%%%%%%%%%%%%
%
%%%%%%%%%%%%%%
\brem{
In Remark~\ref{rem-sum}, we mentioned the relationship
$r_{k,n}^++r_{k,n}^-+r_{k,n}^0=r_{k,n}.$
This one can be checked by adding~\refp{law-bridge},
\refp{law-positive} and~\refp{law-negative} after noticing that
$$
r_{k,n}^0=2\ind_{\cE}(k)\ind_{\cE}(n)a_n(pq)^{n/2}
=\ind_{\cE}(n-k)\ind_{\cE}(n)2p\,a_n(pq)^{n/2}
+\ind_{\cE}(k)\ind_{\cE}(n)2q\,a_n(pq)^{n/2},
$$
the foregoing sum coincides with~\refp{lawfree}.
}
%%%%%%%%%%%%%%

%%%%%%%%%%%%%%%%%%%%%%%%%%%%%%%%%%%%%%%%%%%%%%%%%%%%%%%%%%%%%%%%%%%%%%%%%%%%
%%%%%%%%%%%%%%%%%%%%%%%%%%%%%%%%%%%%%%%%%%%%%%%%%%%%%%%%%%%%%%%%%%%%%%%%%%%%
%%%%%%%%%%%%%%%%%%%%%%%%%%%%%%%%%%%%%%%%%%%%%%%%%%%%%%%%%%%%%%%%%%%%%%%%%%%%
\section{Examples}

In this section, we provide explicit values for $r_{k,n},r_{k,n}^0,r_{k,n}^+$
for $k\le n\le 8$. They can be deduced from our general results. We propose to
compute them by enumerating all the possible paths of the random walk.
We recall some properties which allow to check our results and which prevent us
from considering certain cases. Knowing $r_{k,n}$, we find $r_{n-k,n}$ by
interchanging $p$ and $q$; hence, we only need the values of $r_{k,n}$
for $k\le n/2$. For $k=0$ and $k=1$, we have
$r_{0,n}=1-2p\sum_{i\in\cE:\atop i\le n-1} a_i (pq)^{i/2}$ and
$r_{1,n}= 2 a_{n-1}\, p^{(n+1)/2}q^{(n-1)/2}$. For even $n$ and odd $k$,
$r_{k,n}=0$. For even integers $n,k$ such that $0\le k\le n-1$,
$r_{k-1,n-1}+r_{k,n-1}=r_{k,n}$. The numbers $r_{k,n}^0$ do not vanish only
for even integers $k$ and $n$.
We also have $r_{0,n}^+=r_{1,n}^+=0$, $r_{2,n}^+=p^2r_{0,n-2}^0$ and
if $n-k$ is odd, then $r_{k,n}^+=0$.
We give all the details in the case of the unconditioned random
walk, and only the results for the conditioned cases.

%%%%%%%%%%%%%%%%%%%%%%%%%%%%%%%%%%%%%%%%%%%%%%%%%%%%%%%%%%%%%%%%%%%%%%%%%%%%
\subsection{Case $n=1$}

For the unconditioned case, we simply have
$$
r_{0,1}=\P\{S_1\le 0\}=\P\{S_1=-1\}=q.
$$
The probability distribution of $T_1$ is given by $r_{0,1}=q,$ $r_{1,1}=p$.
Concerning the conditioned case, we have  $r_{0,1}^+=0,$ $r_{1,1}^+=p$.

%%%%%%%%%%%%%%%%%%%%%%%%%%%%%%%%%%%%%%%%%%%%%%%%%%%%%%%%%%%%%%%%%%%%%%%%%%%%
\subsection{Case $n=2$}

For the unconditioned case, we simply have
$$
r_{0,2}=\P\{S_1\le 0,S_2\le 0\}=\P\{S_1=-1\}=q.
$$
The probability distribution of $T_2$ is given by
$$
r_{0,2}=q,\quad r_{1,2}=0,\quad r_{2,2}=p.
$$
We can check that $r_{0,2}=r_{0,1}$. Concerning the conditioned cases, we have
$$
r_{0,2}^0=r_{2,2}^0=pq,\quad r_{0,2}^+=r_{1,2}^+=0,\quad r_{2,2}^+=p^2.
$$

%%%%%%%%%%%%%%%%%%%%%%%%%%%%%%%%%%%%%%%%%%%%%%%%%%%%%%%%%%%%%%%%%%%%%%%%%%%%
\subsection{Case $n=3$}

For the unconditioned case, we have
\begin{align*}
r_{0,3}
&=
\P\{S_1\le 0,S_2\le 0,S_3\le 0\}
=\P\{S_1=-1,S_2=0,S_3=-1\}+\P\{S_1=-1,S_2=-2\}
\\
&=
pq^2+q^2=q^2(p+1),
\\
r_{1,3}
&=
\P\{S_1=-1,S_2=0,S_3=1\}=p^2q.
\end{align*}
The probability distribution of $T_3$ is given by
$$
r_{0,3}=q^2(p+1),\quad r_{1,3}=p^2q,\quad r_{2,3}=pq^2,\quad r_{3,3}=p^2(q+1).
$$
As a check, \refp{lawfree-part-cases} yields
$$
r_{0,3}=1-2p[a_0+a_2pq]=1-p-p^2q=q(1-p^2)=q^2(p+1)
$$
which confirms the result on $r_{0,3}$ for instance.
Concerning the conditioned case, we have
$$
r_{0,3}^+=r_{2,3}^+=0,\quad r_{1,3}^+=p^2q,\quad r_{3,3}^+=p^2(q+1).
$$

%%%%%%%%%%%%%%%%%%%%%%%%%%%%%%%%%%%%%%%%%%%%%%%%%%%%%%%%%%%%%%%%%%%%%%%%%%%%
\subsection{Case $n=4$}

For the unconditioned case, we have
\begin{align*}
r_{0,4}
&=
\P\{S_1\le 0,S_2\le 0,S_3\le 0,S_4\le 0\}
\\
&=
\P\{S_1=-1,S_2=0,S_3=-1\}+\P\{S_1=-1,S_2=-2\}
\\
&=
pq^2+q^2=q^2(p+1),
\\
r_{2,4}
&=
\P\{S_1=1,S_2=0,S_3=-1\}+\P\{S_1=-1,S_2=0,S_3=1\}
\\
&=
pq^2+p^2q=pq.
\end{align*}
The probability distribution of $T_4$ is given by
$$
r_{0,4}=q^2(p+1),\quad r_{1,4}=0,\quad r_{2,4}=pq,\quad r_{3,4}=0,\quad
r_{4,4}=p^2(q+1).
$$
We can check that $r_{0,4}=r_{0,3}$ and $r_{2,4}=r_{2,2}r_{0,2}=r_{1,3}+r_{2,3}.$
Concerning the conditioned cases, we have
$$
r_{0,4}^0=r_{2,4}^0=r_{4,4}^0=2p^2q^2,\quad r_{0,4}^+=r_{1,4}^+=r_{3,4}^+=0,
\quad r_{2,4}^+=p^3q,\quad r_{4,4}^+=p^3(2q+1).
$$
%%%%%%%%%%%%%%%%%%%%%%%%%%%%%%%%%%%%%%%%%%%%%%%%%%%%%%%%%%%%%%%%%%%%%%%%%%%%
\subsection{Case $n=5$}

For the unconditioned case, we have
\begin{align*}
r_{0,5}
&=
\P\{S_1\le 0,S_2\le 0,S_3\le 0,S_4\le 0,S_5\le 0\}
\\
&=
\P\{S_1=-1,S_2=0,S_3=-1,S_4=0,S_5=-1\}
\\
&\hphantom{=\,}
+\P\{S_1=-1,S_2=0,S_3=-1,S_4=-2\}
\\
&\hphantom{=\,}
+\P\{S_1=-1,S_2=-2,S_3=-1,S_4=0,S_5=-1\}
\\
&\hphantom{=\,}
+\P\{S_1=-1,S_2=-2,S_3=-3\}
\\
&=
q^3(2p^2+2p+1),
\\[-2ex]
\lqn{}
r_{1,5}
&=
\P\{S_1=-1,S_2=0,S_3=-1,S_4=0,S_5=1\}
\\
&\hphantom{=\,}
+\P\{S_1=-1,S_2=-2,S_3=-1,S_4=0,S_5=1\}
\\
&=
2p^3q^2,
\\[-2ex]
\lqn{}
r_{2,5}
&=
\P\{S_1=-1,S_2=0,S_3=1,S_4=0,S_5=-1\}
\\
&\hphantom{=\,}
+\P\{S_1=1,S_2=0,S_3=-1,S_4=0,S_5=-1\}
\\
&\hphantom{=\,}
+\P\{S_1=1,S_2=0,S_3=-1,S_4=-2\}
\\
&=
2p^2q^3+pq^3=pq^3(2p+1).
\end{align*}
The probability distribution of $T_5$ is given by
\begin{gather*}
r_{0,5}=q^3(2p^2+2p+1),\quad
r_{1,5}=2p^3q^2,\quad
r_{2,5}=pq^3(2p+1),
\\
r_{3,5}=p^3q(2q+1),\quad
r_{4,5}=2p^2q^3,\quad
r_{5,5}=p^3(2q^2+2q+1).
\end{gather*}
As a check, \refp{lawfree-part-cases} yields
\begin{align*}
r_{0,5}
&=
1-2p[a_0+a_2pq+a_4(pq)^2]=r_{0,3}-2p^3q^2
\\
&=
q^2(p+1-2p^3)
=q^2(1-p)(1+2p+2p^2)=q^3(2p^2+2p+1)
\end{align*}
which confirms the result on $r_{0,5}$.
On the other hand,~\refp{lawfree-odd-case} yields
\begin{align*}
r_{2,5}
&=
2q[a_2pq+a_4(pq)^2]-4a_0a_2(pq)^2
\\
&=
pq^2+2p^2q^3-p^2q^2
=pq^2(1+2pq-p)=pq^3(2p+1)
\end{align*}
which confirms the result on $r_{2,5}$.
Concerning the conditioned case, we have
$$
r_{0,5}^+=r_{2,5}^+=r_{4,5}^+=0,\quad r_{1,5}^+=2p^3q^2,
\quad r_{3,5}^+=p^3q(2q+1),\quad r_{5,5}^+=p^3(2q^2+2q+1).
$$

%%%%%%%%%%%%%%%%%%%%%%%%%%%%%%%%%%%%%%%%%%%%%%%%%%%%%%%%%%%%%%%%%%%%%%%%%%%%
\subsection{Case $n=6$}

For the unconditioned case, we have
\begin{align*}
r_{0,6}
&=
\P\{S_1\le 0,S_2\le 0,S_3\le 0,S_4\le 0,S_5\le 0,S_6\le 0\}
\\
&=
\P\{S_1=-1,S_2=0,S_3=-1,S_4=0,S_5=-1\}
\\
&\hphantom{=\,}
+\P\{S_1=-1,S_2=0,S_3=-1,S_4=-2\}
\\
&\hphantom{=\,}
+\P\{S_1=-1,S_2=-2,S_3=-1,S_4=0,S_5=-1\}
\\
&\hphantom{=\,}
+\P\{S_1=-1,S_2=-2,S_3=-1,S_4=-2\}
\\
&\hphantom{=\,}
+\P\{S_1=-1,S_2=-2,S_3=-3\}
\\
&=
q^3(2p^2+2p+1),
\\
\lqn{}
r_{2,6}
&=
\P\{S_1=-1,S_2=0,S_3=-1,S_4=0,S_5=1\}
\\
&\hphantom{=\,}
+\P\{S_1=-1,S_2=-2,S_3=-1,S_4=0,S_5=1\}
\\
&\hphantom{=\,}
+\P\{S_1=-1,S_2=0,S_3=1,S_4=0,S_5=-1\}
\\
&\hphantom{=\,}
+\P\{S_1=1,S_2=0,S_3=-1,S_4=0,S_5=-1\}
\\
&\hphantom{=\,}
+\P\{S_1=1,S_2=0,S_3=-1,S_4=-2\}
\\
&=
2p^3q^2+2p^2q^3+pq^3=pq^2(p+1).
\end{align*}
The probability distribution of $T_6$ is given by
\begin{gather*}
r_{0,6}=q^3(2p^2+2p+1),\quad
r_{1,6}=0,\quad
r_{2,6}=pq^2(p+1),\quad
r_{3,6}=0,
\\
r_{4,6}=p^2q(q+1),\quad
r_{5,6}=0,\quad
r_{6,6}=p^3(2q^2+2q+1).
\end{gather*}
We can check that $r_{0,6}=r_{0,5}$ and $r_{2,6}=r_{2,2}r_{0,4}=r_{1,5}+r_{2,5}.$
On the other hand,~\refp{lawfree-even-case} yields
\begin{align*}
r_{2,6}
&=
2p[a_4(pq)^2+a_6(pq)^3]-4a_0a_4(pq)^3+2q[a_2(pq)+a_4(pq)^2+a_6(pq)^3]
\\
&\hphantom{=\,}
-4[a_0a_2(pq)^2+(a_0a_4+a_2^2)(pq)^3]
\\
&=
2a_2pq^2+2[a_4(pq)^2+a_6(pq)^3]-4[a_0a_2(pq)^2+a_2^2(pq)^3]
\\
&=
pq^2+(pq)^2=pq^2(p+1)
\end{align*}
which confirms the result on $r_{2,6}$.
Concerning the conditioned cases, we have
\begin{gather*}
r_{0,6}^0=r_{2,6}^0=r_{4,6}^0=r_{6,6}^0=5p^3q^3,
\\
r_{0,6}^+=r_{1,6}^+=r_{3,6}^+=r_{5,6}^+=0, \quad r_{2,6}^+=2p^4q^2,
\quad r_{4,6}^+=p^4q(3q+1),\quad r_{6,6}^+=p^4(5q^2+3q+1).
\end{gather*}

%%%%%%%%%%%%%%%%%%%%%%%%%%%%%%%%%%%%%%%%%%%%%%%%%%%%%%%%%%%%%%%%%%%%%%%%%%%%
\subsection{Case $n=7$}

For the unconditioned case, we have
\begin{align*}
r_{0,7}
&=
\P\{S_1\le 0,S_2\le 0,S_3\le 0,S_4\le 0,S_5\le 0,S_6\le 0,S_7\le 0\}
\\
&=
\P\{S_1=-1,S_2=0,S_3=-1,S_4=0,S_5=-1,S_6=0,S_7=-1\}
\\
&\hphantom{=\,}
+\P\{S_1=-1,S_2=0,S_3=-1,S_4=0,S_5=-1,S_6=-2\}
\\
&\hphantom{=\,}
+\P\{S_1=-1,S_2=0,S_3=-1,S_4=-2,S_5=-1,S_6=0,S_7=-1\}
\\
&\hphantom{=\,}
+\P\{S_1=-1,S_2=0,S_3=-1,S_4=-2,S_5=-1,S_6=-2\}
\\
&\hphantom{=\,}
+\P\{S_1=-1,S_2=0,S_3=-1,S_4=-2,S_5=-3\}
\\
&\hphantom{=\,}
+\P\{S_1=-1,S_2=-2,S_3=-1,S_4=0,S_5=-1,S_6=0,S_7=-1\}
\\
&\hphantom{=\,}
+\P\{S_1=-1,S_2=-2,S_3=-1,S_4=0,S_5=-1,S_6=-2\}
\\
&\hphantom{=\,}
+\P\{S_1=-1,S_2=-2,S_3=-1,S_4=-2,S_5=-1,S_6=0,S_7=-1\}
\\
&\hphantom{=\,}
+\P\{S_1=-1,S_2=-2,S_3=-1,S_4=-2,S_5=-1,S_6=-2\}
\\
&\hphantom{=\,}
+\P\{S_1=-1,S_2=-2,S_3=-1,S_4=-2,S_5=-3\}
\\
&\hphantom{=\,}
+\P\{S_1=-1,S_2=-2,S_3=-3,S_4=-2,S_5=-1,S_6=0,S_7=-1\}
\\
&\hphantom{=\,}
+\P\{S_1=-1,S_2=-2,S_3=-3,S_4=-2,S_5=-1,S_6=-2\}
\\
&\hphantom{=\,}
+\P\{S_1=-1,S_2=-2,S_3=-3,S_4=-2,S_5=-3\}
\\
&\hphantom{=\,}
+\P\{S_1=-1,S_2=-2,S_3=-3,S_4=-4\}
\\
&=
q^4(5p^3+5p^2+3p+1)
\\ \lqn{}
r_{1,7}
&=
\P\{S_1=-1,S_2=0,S_3=-1,S_4=0,S_5=-1,S_6=0,S_7=1\}
\\
&\hphantom{=\,}
+\P\{S_1=-1,S_2=0,S_3=-1,S_4=-2,S_5=-1,S_6=0,S_7=1\}
\\
&\hphantom{=\,}
+\P\{S_1=-1,S_2=-2,S_3=-1,S_4=0,S_5=-1,S_6=0,S_7=1\}
\\
&\hphantom{=\,}
+\P\{S_1=-1,S_2=0,S_3=-1,S_4=-2,S_5=-1,S_6=0,S_7=1\}
\\
&\hphantom{=\,}
+\P\{S_1=-1,S_2=-2,S_3=-3,S_4=-2,S_5=-1,S_6=0,S_7=1\}
\\
&=
5p^4q^3
\\ \lqn{}
r_{2,7}
&=
\P\{S_1=-1,S_2=0,S_3=-1,S_4=0,S_5=1,S_6=0,S_7=-1\}
\\
&\hphantom{=\,}
+\P\{S_1=-1,S_2=-2,S_3=-1,S_4=0,S_5=1,S_6=0,S_7=-1\}
\\
&\hphantom{=\,}
+\P\{S_1=-1,S_2=0,S_3=1,S_4=0,S_5=-1,S_6=0,S_7=-1\}
\\
&\hphantom{=\,}
+\P\{S_1=-1,S_2=0,S_3=1,S_4=0,S_5=-1,S_6=-2\}
\\
&\hphantom{=\,}
+\P\{S_1=1,S_2=0,S_3=-1,S_4=0,S_5=-1,S_6=0,S_7=-1\}
\\
&\hphantom{=\,}
+\P\{S_1=1,S_2=0,S_3=-1,S_4=0,S_5=-1,S_6=-2\}
\\
&\hphantom{=\,}
+\P\{S_1=1,S_2=0,S_3=-1,S_4=-2,S_5=-1,S_6=0,S_7=-1\}
\\
&\hphantom{=\,}
+\P\{S_1=1,S_2=0,S_3=-1,S_4=-2,S_5=-1,S_6=-2\}
\\
&\hphantom{=\,}
+\P\{S_1=1,S_2=0,S_3=-1,S_4=-2,S_5=-3\}
\\
&=
pq^4(5p^2+3p+1)
\\
\lqn{}
r_{3,7}
&=
\P\{S_1=-1,S_2=0,S_3=-1,S_4=0,S_5=1,S_6=0,S_7=1\}
\\
&\hphantom{=\,}
+\P\{S_1=-1,S_2=-2,S_3=-1,S_4=0,S_5=1,S_6=0,S_7=1\}
\\
&\hphantom{=\,}
+\P\{S_1=-1,S_2=0,S_3=1,S_4=0,S_5=-1,S_6=0,S_7=1\}
\\
&\hphantom{=\,}
+\P\{S_1=1,S_2=0,S_3=-1,S_4=0,S_5=-1,S_6=0,S_7=1\}
\\
&\hphantom{=\,}
+\P\{S_1=1,S_2=0,S_3=-1,S_4=-2,S_5=-1,S_6=0,S_7=1\}
\\
&\hphantom{=\,}
+\P\{S_1=-1,S_2=0,S_3=-1,S_4=0,S_5=1,S_6=2\}
\\
&\hphantom{=\,}
+\P\{S_1=-1,S_2=-2,S_3=-1,S_4=0,S_5=1,S_6=2\}
\\
&=
p^4q^2(5q+2)
\end{align*}
The probability distribution of $T_7$ is given by
\begin{gather*}
r_{0,7}=q^4(5p^3+5p^2+3p+1),\quad
r_{1,7}=5p^4q^3,\quad
r_{2,7}=pq^4(5p^2+3p+1),\quad
r_{3,7}=p^4q^2(5q+2),
\\
r_{4,7}=p^2q^4(5p+2),\quad
r_{5,7}=p^4q(5q^2+3q+1),\quad
r_{6,7}=5p^3q^4,\quad
r_{7,7}=p^4(5q^3+5q^2+3q+1).
\end{gather*}
As a check, \refp{lawfree-part-cases} yields
\begin{align*}
r_{0,7}
&=
1-2p[a_0+a_2pq+a_4(pq)^2+a_6(pq)^3]
\\
&=
r_{0,5}-5p^4q^3=q^3(2p^2+2p+1-5p^4)
\\
&=
q^3(1-p)(5p^3+5p^2+3p+1)=q^4(5p^3+5p^2+3p+1)
\end{align*}
which confirms the result on $r_{0,7}$.
On the other hand,~\refp{lawfree-odd-case} yields
\begin{align*}
r_{2,7}
&=
2q[a_2pq+a_4(pq)^2+a_6(pq)^3]-4[a_0a_2(pq)^2+(a_0a_4+a_2^2)(pq)^3]
\\
&=
r_{2,5}+5p^3q^4-3p^3q^3=pq^3(2p+1+5p^2q-3p^2)
\\
&=
pq^3((1-p)(1+3p)+5p^2q)=pq^4(5p^2+3p+1)
\\
r_{3,7}
&=
2p[a_4(pq)^2+a_6(pq)^3]-4a_0a_4(pq)^3
\\
&=
2p^3q^2+5p^4q^3-2p^3q^3=2p^4q^2+5p^4q^3=p^4q^2(5q+2)
\end{align*}
which confirms the results on $r_{2,7}$ and $r_{3,7}$.
Concerning the conditioned case, we have
\begin{gather*}
r_{0,7}^+=r_{2,7}^+=r_{4,7}^+=r_{6,7}^+=0,\quad r_{1,7}^+=5p^4q^3,
\quad r_{3,7}^+=p^4q^2(5q+2),
\\
r_{5,7}^+=p^4q(5q^2+3q+1), \quad r_{7,7}^+=p^4(5q^3+5q^2+3q+1).
\end{gather*}

%%%%%%%%%%%%%%%%%%%%%%%%%%%%%%%%%%%%%%%%%%%%%%%%%%%%%%%%%%%%%%%%%%%%%%%%%%%%
\subsection{Case $n=8$}

For the unconditioned case, we have
\begin{align*}
r_{0,8}&=\P\{S_1\le 0,S_2\le 0,S_3\le 0,S_4\le 0,S_5\le 0,S_6\le 0,S_7\le 0,S_8\le 0\}
\\
&=
\P\{S_1=-1,S_2=0,S_3=-1,S_4=0,S_5=-1,S_6=0,S_7=-1\}
\\
&\hphantom{=\,}
+\P\{S_1=-1,S_2=0,S_3=-1,S_4=0,S_5=-1,S_6=-2\}
\\
&\hphantom{=\,}
+\P\{S_1=-1,S_2=0,S_3=-1,S_4=-2,S_5=-1,S_6=0,S_7=-1\}
\\
&\hphantom{=\,}
+\P\{S_1=-1,S_2=0,S_3=-1,S_4=-2,S_5=-1,S_6=-2\}
\\
&\hphantom{=\,}
+\P\{S_1=-1,S_2=0,S_3=-1,S_4=-2,S_5=-3\}
\\
&\hphantom{=\,}
+\P\{S_1=-1,S_2=-2,S_3=-1,S_4=0,S_5=-1,S_6=0,S_7=-1\}
\\
&\hphantom{=\,}
+\P\{S_1=-1,S_2=-2,S_3=-1,S_4=0,S_5=-1,S_6=-2\}
\\
&\hphantom{=\,}
+\P\{S_1=-1,S_2=-2,S_3=-1,S_4=-2,S_5=-1,S_6=0,S_7=-1\}
\\
&\hphantom{=\,}
+\P\{S_1=-1,S_2=-2,S_3=-1,S_4=-2,S_5=-3\}
\\
&\hphantom{=\,}
+\P\{S_1=-1,S_2=-2,S_3=-1,S_4=-2,S_5=-1,S_6=-2\}
\\
&\hphantom{=\,}
+\P\{S_1=-1,S_2=-2,S_3=-3,S_4=-2,S_5=-1,S_6=0,S_7=-1\}
\\
&\hphantom{=\,}
+\P\{S_1=-1,S_2=-2,S_3=-3,S_4=-2,S_5=-1,S_6=-2\}
\\
&\hphantom{=\,}
+\P\{S_1=-1,S_2=-2,S_3=-3,S_4=-2,S_5=-3\}
\\
&\hphantom{=\,}
+\P\{S_1=-1,S_2=-2,S_3=-3,S_4=-4\}
\\
&=
q^4(5p^3+5p^2+3p+1)
\\ \lqn{}
r_{2,8}&= \P\{S_1=-1,S_2=0,S_3=-1,S_4=0,S_5=-1,S_6=0,S_7=1\}
\\
&\hphantom{=\,}
+\P\{S_1=-1,S_2=0,S_3=-1,S_4=-2,S_5=-1,S_6=0,S_7=1\}
\\
&\hphantom{=\,}
+\P\{S_1=-1,S_2=-2,S_3=-1,S_4=0,S_5=-1,S_6=0,S_7=1\}
\\
&\hphantom{=\,}
+\P\{S_1=-1,S_2=-2,S_3=-1,S_4=-2,S_5=-1,S_6=0,S_7=1\}
\\
&\hphantom{=\,}
+\P\{S_1=-1,S_2=-2,S_3=-3,S_4=-2,S_5=-1,S_6=0,S_7=1\}
\\
&\hphantom{=\,}
+\P\{S_1=-1,S_2=0,S_3=-1,S_4=0,S_5=1,S_6=0,S_7=-1\}
\\
&\hphantom{=\,}
+\P\{S_1=-1,S_2=-2,S_3=-1,S_4=0,S_5=1,S_6=0,S_7=-1\}
\\
&\hphantom{=\,}
+\P\{S_1=-1,S_2=0,S_3=1,S_4=0,S_5=-1,S_6=0,S_7=-1\}
\\
&\hphantom{=\,}
+\P\{S_1=-1,S_2=0,S_3=1,S_4=0,S_5=-1,S_6=-2\}
\\
&\hphantom{=\,}
+\P\{S_1=1,S_2=0,S_3=-1,S_4=0,S_5=-1,S_6=0,S_7=-1\}
\\
&\hphantom{=\,}
+\P\{S_1=1,S_2=0,S_3=-1,S_4=0,S_5=-1,S_6=-2\}
\\
&\hphantom{=\,}
+\P\{S_1=1,S_2=0,S_3=-1,S_4=-2,S_5=-1,S_6=0,S_7=-1\}
\\
&\hphantom{=\,}
+\P\{S_1=1,S_2=0,S_3=-1,S_4=-2,S_5=-1,S_6=-2\}
\\
&\hphantom{=\,}
+\P\{S_1=1,S_2=0,S_3=-1,S_4=-2,S_5=-3\}
\\
&=
5p^4q^3+5p^3q^4+3p^2q^4+pq^4=pq^3[5p^3+(5p^2+3p+1)(1-p)]
\\
&=
pq^3(2p^2+2p+1).
\\ \lqn{}
r_{4,8}&= \P\{S_1=-1,S_2=0,S_3=-1,S_4=0,S_5=1,S_6=0,S_7=1\}
\\
&\hphantom{=\,}
+\P\{S_1=-1,S_2=-2,S_3=-1,S_4=0,S_5=1,S_6=0,S_7=1\}
\\
&\hphantom{=\,}
+\P\{S_1=-1,S_2=0,S_3=-1,S_4=0,S_5=1,S_6=2\}
\\
&\hphantom{=\,}
+\P\{S_1=-1,S_2=-2,S_3=-1,S_4=0,S_5=1,S_6=2\}
\\
&\hphantom{=\,}
+\P\{S_1=-1,S_2=0,S_3=1,S_4=0,S_5=-1,S_6=0,S_7=1\}
\\
&\hphantom{=\,}
+\P\{S_1=-1,S_2=0,S_3=1,S_4=0,S_5=1,S_6=0,S_7=-1\}
\\
&\hphantom{=\,}
+\P\{S_1=-1,S_2=0,S_3=1,S_4=2,S_5=1,S_6=0,S_7=-1\}
\\
&\hphantom{=\,}
+\P\{S_1=1,S_2=0,S_3=-1,S_4=0,S_5=-1,S_6=0,S_7=1\}
\\
&\hphantom{=\,}
+\P\{S_1=1,S_2=0,S_3=-1,S_4=-2,S_5=-1,S_6=0,S_7=1\}
\\
&\hphantom{=\,}
+\P\{S_1=1,S_2=0,S_3=-1,S_4=0,S_5=1,S_6=0,S_7=-1\}
\\
&\hphantom{=\,}
+\P\{S_1=1,S_2=0,S_3=1,S_4=0,S_5=-1,S_6=0,S_7=-1\}
\\
&\hphantom{=\,}
+\P\{S_1=1,S_2=0,S_3=1,S_4=0,S_5=-1,S_6=-2\}
\\
&\hphantom{=\,}
+\P\{S_1=1,S_2=2,S_3=1,S_4=0,S_5=-1,S_6=0,S_7=-1\}
\\
&\hphantom{=\,}
+\P\{S_1=1,S_2=2,S_3=1,S_4=0,S_5=-1,S_6=-2\}
\\
&=
5p^4q^3+5p^3q^4+2p^4q^2+2p^2q^4=5p^3q^3(p+q)+2p^2q^2(p^2+q^2)
\\
&=
p^2q^2(-p^2+p+2).
\end{align*}
The probability distribution of $T_8$ is given by
\begin{gather*}
r_{0,8}=q^4(5p^3+5p^2+3p+1),\quad
r_{1,8}=0,\quad
r_{2,8}=pq^3(2p^2+2p+1),
\\
r_{3,8}=0,\quad
r_{4,8}=p^2q^2(-p^2+p+2),\quad
r_{5,8}=0,
\\
r_{6,8}=p^3q(2q^2+2q+1),\quad
r_{7,8}=0,\quad
r_{8,8}=p^4(5q^3+5q^2+3q+1).
\end{gather*}
We can check that $r_{0,8}=r_{0,7}$, $r_{2,8}=r_{2,2}r_{0,6}=r_{1,7}+r_{2,7}$,
and $r_{4,8}=r_{4,4}r_{0,4}=r_{3,7}+r_{4,7}$.
On the other hand,~\refp{lawfree-even-case} yields
\begin{align*}
r_{2,8}
&=
2p[a_6(pq)^3+a_8(pq)^4]-4a_0a_6(pq)^4+2q[a_2(pq)+a_4(pq)^2+a_6(pq)^3+a_8(pq)^4]
\\
&\hphantom{=\,}
-4[a_0a_2(pq)^2+(a_0a_4+a_2^2)(pq)^3+(a_0a_6+2a_2a_4)(pq)^4]
\\
&=
pq^2+2p^2q^3+5(pq)^3+14(pq)^4-(pq)^2-3(pq)^3-14(pq)^4
\\
&
=pq^2+2p^2q^3-p^2q^2+2p^3q^3=pq^3+2p^2q^3+2p^3q^3
=pq^3(2p^2+2p+1)
\\
r_{4,8}
&=
2[a_4(pq)^2+a_6(pq)^3+a_8(pq)^4]-8[a_0a_4(pq)^3+(a_0a_6+a_2a_4)(pq)^4]
\\
&=
2(pq)^2+5(pq)^3+14(pq)^4-4(pq)^3-14(pq)^4
=p^2q^2(pq+2)
\end{align*}
which confirms the results on $r_{2,8}$ and $r_{4,8}$.
Concerning the conditioned cases, we have
\begin{gather*}
r_{0,8}^0=r_{2,8}^0=r_{4,8}^0=r_{6,8}^0=r_{8,8}^0=14p^4q^4,
\\
r_{0,8}^+=r_{1,8}^+=r_{3,8}^+=r_{5,8}^+=r_{7,8}^+=0, \quad r_{2,8}^+=5p^5q^3,
\quad r_{4,8}^+=p^5q^2(7q+2),
\\
r_{6,8}^+=p^5q(9q^2+4q+1),\quad r_{8,8}^+=p^5(14q^3+9q^2+4q+1).
\end{gather*}

%%%%%%%%%%%%%%%%%%%%%%%%%%%%%%%%%%%%%%%%%%%%%%%%%%%%%%%%%%%%%%%%%%%%%%%%%%%%%%%
%%%%%%%%%%%%%%%%%%%%%%%%%       Asymptotics    %%%%%%%%%%%%%%%%%%%%%%%%%%%%%%%%
%%%%%%%%%%%%%%%%%%%%%%%%%%%%%%%%%%%%%%%%%%%%%%%%%%%%%%%%%%%%%%%%%%%%%%%%%%%%%%%
\section{Asymptotics: Brownian motion with a linear drift}

In this part, our aim is to retrieve certain probability distributions
related to the sojourn time in $(0,+\infty)$ of Brownian motion
with a linear drift.

%%%%%%%%%%%%%%%%%%%%%%%%%%%%%%%%%%%%%%%%%%%%%%%%%%%%%%%%%%%%%%%%%%%%%%%%%%%%%%%
%%%%%%%%%%%%%%%%%%%%%%%%%       Rescaling      %%%%%%%%%%%%%%%%%%%%%%%%%%%%%%%%
%%%%%%%%%%%%%%%%%%%%%%%%%%%%%%%%%%%%%%%%%%%%%%%%%%%%%%%%%%%%%%%%%%%%%%%%%%%%%%%
\subsection{Rescaled random walk}

We consider a sequence of random walks $(S_k^\sN)_{k\in\N}$ indexed by $N\in\N$,
defined by
$$
S_k^\sN=S_0^\sN+\sum_{j=1}^k X_j^\sN, \;k\ge 1,
$$
where for each $N\in\N^*$, $(X_j^\sN)_{j\in\N^*}$ is a sequence of independent
Bernoulli variables with jump probabilities depending on $N$ as follows:
$$
\pN=\P\{X_k^\sN=+1\}=\frac12+\frac{\r}{2\sqrt N}
\quad\mbox{and}\quad \qN=\P\{X_k^\sN=-1\}=\frac12-\frac{\r}{2\sqrt N},
$$
$\r$ being a fixed parameter. We also define the centered random walk
$(\tS_k^\sN)_{k\in\N}$ as
$$
\tS_k^\sN=S_k^\sN-\E(S_k^\sN)=S_k^\sN-\r\,\frac{k}{\sqrt N}.
$$
Let $T_n^\sN$ be the corresponding sojourn time in $\Z^+$:
$T_n^\sN=\sum_{j=1}^n \d_j^\sN$ with
$$
\d_j^\sN=\left\{\begin{array}{ll}
1 & \mbox{if $(S_j^\sN>0)$ or $(S_j^\sN=0$ and $S_{j-1}^\sN>0)$,}
\\
0 & \mbox{if $(S_j^\sN<0)$ or $(S_j^\sN=0$ and $S_{j-1}^\sN<0)$.}
\end{array}\right.
$$
Let us introduce the rescaled random walks
$$
B_t^\sN=\frac{1}{\sqrt N}\,S_{[Nt]}^\sN \quad\mbox{and}\quad
\tB_t^\sN=B_t^\sN-\E(B_t^\sN)=\frac{1}{\sqrt N}\,\tS_{[Nt]}^\sN
$$
defined on continuous time $t\ge 0$.
We have
$$
B_t^\sN=\tB_t^\sN+\r\,\frac{[Nt]}{N},\; t\ge 0.
$$
Donsker's theorem (see, e.g., \cite[p. 68]{bill}) asserts that the sequence of processes
$(\tB_t^\sN)_{t\ge 0}$, $N\in\N$, weakly converges to the standard linear
Brownian motion $(\tB_t)_{t\ge 0}$ and the sequence of processes
$(\tB_t^\sN)_{t\ge 0}$, $N\in\N$, weakly converges to the drifted Brownian
motion $(B_t)_{t\ge 0}$ defined as
$$
B_t=\tB_t+\r t,\; t\ge 0.
$$

We now introduce the sojourn times in $\R^+$ of the processes
$(B_t^\sN)_{t\ge 0}$ and $(B_t)_{t\ge 0}$:
$$
\bT_t^\sN=\int_0^t \ind_{\R^+}(B_s^\sN)\,ds \quad\mbox{and}\quad
\bT_t=\int_0^t \ind_{\R^+}(B_s)\,ds.
$$
We have
\begin{align}
\bT_t^\sN
&=
\sum_{j=0}^{[Nt]} \int_{j/N}^{(j+1)/N} \ind_{\R^+} (S_j^\sN) \,ds
-\int_t^{([Nt]+1)/N} \ind_{\R^+} (S_{[Nt]}^\sN) \,ds
\nonumber\\
&=
\frac1N \sum_{j=0}^{[Nt]}\ind_{\R^+} (S_j^\sN)-\Big(\frac{[Nt]+1}{N}-t\Big)
\ind_{\R^+} (S_{[Nt]}^\sN)
\nonumber\\
&=
\frac1N \,T_{[Nt]}^\sN +\frac1N \sum_{j=1}^{[Nt]}\Big(\ind_{\R^+}
(S_j^\sN)-\d_j^\sN\Big) +\frac1N \,\ind_{\R^+} (S_0^\sN)
-\Big(\frac{[Nt]+1}{N}-t\Big)\ind_{\R^+} (S_{[Nt]}^\sN).
\label{sojourn-inter}
\end{align}
We compare the sojourn time of the random walk $(S_k^\sN)_{k\in\N}$
and that of the Brownian motion $(B_t)_{t\ge 0}$:
$$
\frac1N \,T_{[Nt]}^\sN -\bT_t
=\left(\frac1N \,T_{[Nt]}^\sN-\bT_t^\sN\right)+\left(\bT_t^\sN-\bT_t\right).
$$
On one hand, for $j\ge 1$,
$$
\ind_{\R^+} (S_j^\sN)-\d_j^\sN=\left\{\begin{array}{ll}
0 & \mbox{if $(S_j^\sN\ne 0)$ or $(S_j^\sN=0$ and $S_{j-1}^\sN>0)$,}
\\
1 & \mbox{if $S_j^\sN=0$ and $S_{j-1}^\sN<0$,}
\end{array}\right.
$$
and then,
$$
|\ind_{\R^+} (S_j^\sN)-\d_j^\sN|\le \ind_{\{0\}}(S_j).
$$
On the other hand, by~\refp{sojourn-inter}, the following estimate holds:
$$
\left|\frac1N \,T_{[Nt]}^\sN-\bT_t^\sN\right|
\le\frac 1N\sum_{j=1}^{[Nt]}\ind_{\{0\}} (S_j^\sN)+\frac 2N.
$$
We observe, by~\refp{elem}, that
$$
\E\bigg[\sum_{j=1}^{[Nt]}\ind_{\{0\}} (S_j^\sN)\bigg]\le
\sum_{j=1}^{\infty} \P\{S_j^\sN=0\}
\le \sum_{j\in\cE} b_j(\pN\qN)^{j/2}
=\frac{1}{A(\pN\qN)}=\frac{\sqrt N}{\rho}=o(N)
$$
which implies that $\frac1N \,T_{[Nt]}^\sN -\bT_t$ tends to 0 in mean.
Now, by Donsker's theorem, the sequence $(\bT_t^\sN)_{N\in\N^*}$ weakly
converges to $\bT_t$ (see~\cite[p. 72]{bill}).
This discussion shows that $\frac1N \,T_{[Nt]}^\sN -\bT_t$ converges to
0 in mean. As a result, the probability distribution of
$\frac1N \,T_{[Nt]}^\sN$ converges to that of $\bT_t$.
In the following subsection, we compute this limit.

%%%%%%%%%%%%%%%%%%%%%%%%%%%%%%%%%%%%%%%%%%%%%%%%%%%%%%%%%%%%%%%%%%%%%%%%%%%%%%%
%%%%%%%%%%%%%%%%%%%%%%%%%        Limiting      %%%%%%%%%%%%%%%%%%%%%%%%%%%%%%%%
%%%%%%%%%%%%%%%%%%%%%%%%%%%%%%%%%%%%%%%%%%%%%%%%%%%%%%%%%%%%%%%%%%%%%%%%%%%%%%%
\subsection{Limiting distribution of the sojourn time}

%
%%%%%%%%%%%%%%
\bth{\label{limit}
The probability distribution function of the sojourn time $\bT_t$ admits
the following expression: for $0<s<t$,
$$
\P\{\bT_t\in ds\}/ds=\left\{\begin{array}{ll}
\dis \frac{1}{\sqrt{2\pi}} \bigg[\r+\frac{1}{2\sqrt{2\pi}}
\int_s^{\infty} \frac{e^{-\r^2z/2}}{z^{3/2}}\,dz\bigg]
\int_{t-s}^{\infty} \frac{e^{-\r^2z/2}}{z^{3/2}}\,dz & \mbox{if $\r\ge 0$,}
\\[2ex]
\dis \frac{1}{\sqrt{2\pi}} \bigg[|\r|+\frac{1}{2\sqrt{2\pi}}
\int_{t-s}^{\infty} \frac{e^{-\r^2z/2}}{z^{3/2}}\,dz\bigg]
\int_s^{\infty} \frac{e^{-\r^2z/2}}{z^{3/2}}\,dz & \mbox{if $\r\le 0$.}
\end{array}\right.
$$
}
%%%%%%%%%%%%%%
%
The integral $\int_{\s}^{\infty} \frac{e^{-\r^2z/2}}{z^{3/2}}\,dz$ can be expressed
by means of the error function according as
$$
\int_{\s}^{\infty} \frac{e^{-\r^2z/2}}{z^{3/2}}\,dz
=\frac{2}{\sqrt{\s}}\,e^{-\r^2\s/2}-\sqrt{2\pi}\,\r\,\mathrm{Erfc}(\r\sqrt{\s/2}).
$$
Let us point out that the density of $\bT_t$ is known under another form,
see formula 2.1.4.8 in~\cite{borodin}.

\dem
We shall assume that $\r\ge 0$, the case $\r\le 0$ being quite similar.
We begin by computing the limit of the following probability as $N\to \infty$:
$$
\P\{s_1< \frac 1N\,T_{[Nt]}^\sN\le s_2\}=\P\{[Ns_1]< T_{[Nt]}^\sN\le[Ns_2]\}
=\sum_{k=[Ns_1]+1}^{[Ns_2]} \P\{T_{[Nt]}^\sN=k\}.
$$
%%%%%%%%%%%%%%%%%%%%%%%%%%%%%%%%%%%%%%%%%%%%%%%%%%%%%%%%%%%%%%%%%%%%%%%%%%%%%%%
\textsl{$\bullet$ Case where $[Nt]$ is even.}
Suppose that $[Nt]$ is even and  set $\a_i^\sN= a_i(\pN\qN)^{i/2}$
with $\pN\qN=\frac14(1-\r^2/N)$. By using~\refp{lawfree-part-cases}
and~\refp{law-product}, for even integer $k$ such that $[Ns_1]< k\le[Ns_2]$, we get
\begin{align}
\P\{T_{[Nt]}^\sN=k\}
&=
\P\{T_k^\sN=k\}\P\{T_{[Nt]-k}^\sN=0\}
\nonumber\\
&=
\bigg[1-\qN\sum_{i\in\cE:\atop 0\le i\le k-1} \a_i^\sN\bigg]
\bigg[1-\pN\sum_{i\in\cE:\atop 0\le i\le [Nt]-k-1} \a_i^\sN\bigg].
\label{asymp-product1}
\end{align}
Since $\r$ is assumed to be non negative, we have by~\refp{elem}
$$
\sum_{i\in\cE} \a_i^\sN =\frac{1-A(\sqrt{\pN\qN})}{2\pN\qN}
=\frac{1-|\pN-\qN|}{2\pN\qN}=\frac{1}{\pN},
$$
and we can rewrite~\refp{asymp-product1} as
\beq
\P\{T_{[Nt]}^\sN=k\}=\bigg[1-\frac{\qN}{\pN}+
\qN\sum_{i\in\cE:\atop i\ge k} \a_i^\sN\bigg]
\bigg[\pN\sum_{i\in\cE:\atop i\ge [Nt]-k} \a_i^\sN\bigg].
\label{asymp-product2}
\eeq
We aim to evaluate the limit of the foregoing quantity as $N\to\infty$.
For this, we need an asymptotic for the sum
$\sum_{i\in\cE:\atop i\ge k} \a_i^\sN$ as $k,N\to \infty$ with
$c_1N\le k\le c_2N$ for any $c_1,c_2>0$.
%
%%%%%%%%%%%%%%
\blm{\label{lemma}
The following asymptotics holds: for any $c_1,c_2>0$ such that $c_1<c_2$,
\beq
\sum_{i\in\cE:\atop i\ge k} \a_i^\sN
\underset{k,N\to \infty\atop c_1N\le k\le c_2N}{\sim}
\sqrt{\frac{2}{\pi N}} \int_{k/N}^{\infty} \frac{e^{-\r^2z/2}}{z^{3/2}}\,dz.
\label{asymptotic}
\eeq
}
%%%%%%%%%%%%%%
%
The proof of this lemma is postponed to the appendix.
Set, for $z>0$ and $s\in (0,t)$,
$$
\f(z)=\frac{e^{-\r^2z/2}}{z^{3/2}}
\quad\mbox{and}\quad
\psi(s,t)=\sqrt{\frac{2}{\pi}} \bigg[\r+\frac{1}{2\sqrt{2\pi}}\int_s^{\infty}
\f(z)\,dz\bigg] \int_{t-s}^{\infty} \f(z)\,dz.
$$
In the light of~\refp{asymp-product2} and~\refp{asymptotic}, for even $k$ such
that $[Ns_1]<k\le [Ns_2]$, we have
\begin{align*}
\P\{T_{[Nt]}^\sN=k\}
&\underset{N\to \infty\atop [Nt]\in\cE}{\sim}
\bigg[\frac{2\r}{\sqrt N}+\frac{1}{\sqrt{2\pi N}}\int_{k/N}^{\infty} \f(z)\,dz\bigg]
\bigg[\frac{1}{\sqrt{2\pi N}}\int_{([Nt]-k)/N}^{\infty} \f(z)\,dz\bigg]
\\
&\underset{N\to \infty\atop [Nt]\in\cE}{\sim}
\frac 1N \,\psi(k/n).
\end{align*}
Finally,
$$
\P\{s_1< \frac 1N\,T_{[Nt]}^\sN\le s_2\}
\underset{N\to \infty\atop [Nt]\in\cE}{\sim}
\frac 1N \sum_{k\in\cE:\atop [Ns_1]\le k\le [Ns_2]} \psi(k/N)
\underset{N\to \infty\atop [Nt]\in\cE}{\rightarrow}
\frac 12 \int_{s_1}^{s_2}\psi(s,t)\,ds.
$$

%%%%%%%%%%%%%%%%%%%%%%%%%%%%%%%%%%%%%%%%%%%%%%%%%%%%%%%%%%%%%%%%%%%%%%%%%%%%%%%
\textsl{$\bullet$ Case where $[Nt]$ is odd.}

Assume that $[Nt]$ is odd. Then, by invoking~\refp{rec2}, we obtain
\begin{align*}
\P\{s_1< \frac 1N\,T_{[Nt]}^\sN\le s_2\}
&=
\sum_{k\in\cE:\atop [Ns_1]< k\le [Ns_2]} \P\{T_{[Nt]}^\sN=k\}
+\sum_{k\in\cO:\atop [Ns_1]< k\le [Ns_2]} \P\{T_{[Nt]}^\sN=k\}
\\
\lqn{}
&=
\sum_{k\in\cE:\atop [Ns_1]< k\le [Ns_2]} \P\{T_{[Nt]}^\sN=k\}
+\sum_{k\in\cE:\atop [Ns_1]+1< k\le [Ns_2]+1} \P\{T_{[Nt]}^\sN=k-1\}
\\
&=
\sum_{k=[Ns_1]+1}^{[Ns_2]} \P\{T_{[Nt]+1}^\sN=k\}+\eN
\\
&=
\P\{[Ns_1]<T_{[Nt]+1}^\sN\le [Ns_2]\}+\eN
\end{align*}
where $\eN=\ind_{\cO}([Ns_2])\P\{T_{[Nt]}^\sN=[Ns_2]\}
-\ind_{\cO}([Ns_1])\P\{T_{[Nt]}^\sN=[Ns_1]\}.$
We have
\begin{align*}
|\eN|
&\le
\P\{T_{[Nt]}^\sN=[Ns_1]\}+\P\{T_{[Nt]}^\sN=[Ns_2]\}
\\
&\le
\P\{T_{[Nt]+1}^\sN\in\{[Ns_1],[Ns_1]+1\}\}
+\P\{T_{[Nt]+1}^\sN\in\{[Ns_2],[Ns_2]+1\}\}
\\
&\underset{N\to \infty\atop [Nt]\in\cO}{\sim}
\frac1N\left[\psi\Big(\frac{[Ns_1]}{N}\Big)+\psi\Big(\frac{[Ns_1]+1)}{N}\Big)
+\psi\Big(\frac{[Ns_2]}{N}\Big)+\psi\Big(\frac{[Ns_2]+1)}{N}\Big)\right]
\\
&\underset{N\to \infty\atop [Nt]\in\cO}{\sim}
\frac2N \,[\psi(s_1)+\psi(s_1)].
\end{align*}
Since $[Nt]$ is odd, $[Nt]+1$ is even and we can use the foregoing analysis.
Hence,
$$
\P\{s_1\le \frac 1N\,T_{[Nt]}^\sN\le s_2\}
\underset{N\to \infty\atop [Nt]\in\cO}{\sim}
\frac 12 \int_{s_1}^{s_2}\psi(s,t)\,ds.
$$
We know that $\frac 1N\,T_{[Nt]}^\sN\underset{N\to \infty}{\rightarrow} \bT_t$, then
$$
\P\{s_1\le \bT_t \le s_2\}=\frac 12 \int_{s_1}^{s_2}\psi(s,t)\,ds
\quad\mbox{and}\quad \P\{\bT_t\in ds\}/ds=\frac 12\,\psi(s,t)
$$
which ends up the proof of Theorem~\ref{limit}.
\fin

For $\r=0$, we retrieve the famous Paul L\'evy's arcsine law for
standard Brownian motion: for $s\in (0,t)$,
$$
\P\{\bT_t\in ds\}/ds=\frac{1}{\pi\sqrt{s(t-s)}}.
$$

%
%%%%%%%%%%%%%%
\brem{
Let $t$ tend to $\infty$ in the formula of Theorem~\ref{limit}.
The limiting random variable $\bT_{\infty}$ denotes the total sojourn
time in $\R^+$ and its probability distribution is given by
$$
\P\{\bT_{\infty}\in ds\}/ds=\left\{\begin{array}{ll}
0 & \mbox{if $\r\ge 0$,}
\\
\dis \frac{|\rho|}{\sqrt{2\pi}} \int_s^{\infty} \frac{e^{-\r^2z/2}}{z^{3/2}}\,dz
& \mbox{if $\r<0$.}
\end{array}\right.
$$
From this, we deduce $\P\{\bT_{\infty}=\infty\}=1$ if $\rho\ge 0$. If $\rho<0$,
$$
\P\{\bT_{\infty}<\infty\}
=\frac{|\rho|}{\sqrt{2\pi}} \int_0^{\infty} ds\int_s^{\infty} \frac{e^{-\r^2z/2}}{z^{3/2}}\,dz
=\frac{|\rho|}{\sqrt{2\pi}} \int_0^{\infty} \frac{e^{-\r^2z/2}}{\sqrt z}\,dz=1.
$$
This is in good accordance with the effect of the drift near infinity.
}
%%%%%%%%%%%%%%
%

%
%%%%%%%%%%%%%%
\bth{\label{limit-conditioned}
The probability distribution functions of $\bT_t\ind_{(-\infty,0)}(B_t)$
and $\bT_t\ind_{(0,+\infty)}(B_t)$ admit the following expressions: for $0<s<t$,
\begin{align*}
\P\{\bT_t\in ds,B_t<0\}/ds
&=
\frac{\r^-}{\sqrt{2\pi}} \int_s^t \frac{e^{-\r^2z/2}}{z^{3/2}}\,dz
+\frac{1}{4\pi} \int_s^t \frac{e^{-\r^2u/2}}{u^{3/2}}\,du
\int_{t-u}^{\infty} \frac{e^{-\r^2v/2}}{v^{3/2}}\,dv,
\\
\P\{\bT_t\in ds,B_t>0\}/ds
&=
\frac{\r^+}{\sqrt{2\pi}} \int_{t-s}^t \frac{e^{-\r^2z/2}}{z^{3/2}}\,dz
+\frac{1}{4\pi} \int_{t-s}^t \frac{e^{-\r^2u/2}}{u^{3/2}}\,du
\int_{t-u}^{\infty} \frac{e^{-\r^2v/2}}{v^{3/2}}\,dv.
\end{align*}
}
%%%%%%%%%%%%%%
%
\dem
Set $\a_i=a_i(pq)^{i/2}$. We first begin by rewriting $r_{k,n}^-$ as follows.
By~\refp{law-negative}, we have
\begin{align*}
r_{k,n}^-
&=
\ind_{\cE}(k)\bigg[2q\sum_{i\in\cE:\atop k\le i\le n-1} \a_i
-4pq\sum_{i\in\cE:\atop k\le i\le n-2}
\bigg(\sum_{j\in \cE:\atop j\le i-k} \a_j \a_{i-j}\bigg) \bigg].
\end{align*}
The double sum of the foregoing expression can be written as
\begin{align*}
\sum_{i\in\cE:\atop k\le i\le n-2}
\bigg(\sum_{j\in \cE:\atop j\le i-k} \a_j \a_{i-j}\bigg)
&=
\sum_{j\in \cE:\atop j\le n-k-2} \a_j \bigg(\sum_{i\in\cE:\atop j+k\le i\le n-2}\a_{i-j}\bigg)
=\sum_{j\in \cE:\atop j\le n-k-2} \a_j \bigg(\sum_{i\in\cE:\atop k\le i\le n-j-2}\a_i\bigg)
\\
&=
\sum_{i\in\cE:\atop k\le i\le n-2}\a_i \bigg(\sum_{j\in \cE:\atop j\le n-i-2} \a_j\bigg)
=\sum_{i\in\cE:\atop k\le i\le n-2}\a_i \bigg(\sum_{j\in \cE} \a_j-
\sum_{j\in \cE:\atop j\ge n-i-1} \a_j\bigg)
\\
&=
\frac{1}{2(p\vee q)} \sum_{i\in\cE:\atop k\le i\le n-2}\a_i
-\sum_{i\in\cE:\atop k\le i\le n-2}\a_i\bigg(\sum_{j\in \cE:\atop j\ge n-i-1} \a_j\bigg).
\end{align*}
Therefore,
$$
r_{k,n}^-=\ind_{\cE}(k)\bigg[2q\ind_{\cO}(n)\a_{n-1}
+2(q-p)^+\sum_{i\in\cE:\atop k\le i\le n-2} \a_i
+4pq\sum_{i\in\cE:\atop k\le i\le n-2} \a_i
\bigg(\sum_{j\in \cE:\atop j\ge n-i-1} \a_j\bigg) \bigg].
$$
Recall that $\a_i^\sN=a_i(\pN\qN)^{i/2}$. We now search an asymptotic for the
following probability:
$$
\P\{s_1< \frac 1N\,T_{[Nt]}^\sN\le s_2,B_t^\sN<0\}
=\sum_{k=[Ns_1]+1}^{[Ns_2]} \P\{T_{[Nt]}^\sN=k,S_{[Nt]}<0\}
$$
where
\begin{align*}
\P\{T_{[Nt]}^\sN=k,S_{[Nt]}<0\}
&=
\ind_{\cE}(k)\bigg[2\qN\ind_{\cO}([Nt])\a_{[Nt]-1}^\sN
+2(\qN-\pN)^+\sum_{i\in\cE:\atop k\le i\le [Nt]-2} \a_i^\sN
\\
&\hphantom{=\,}
+4\pN\qN\sum_{i\in\cE:\atop k\le i\le [Nt]-2} \a_i^\sN
\bigg(\sum_{j\in \cE:\atop j\ge [Nt]-i-1} \a_j^\sN\bigg) \bigg].
\end{align*}
Observing that, as $N\to \infty$,
$$
(\qN-\pN)^+\sim \frac{2\rho^-}{\sqrt N},\quad
4\pN\qN \sim 1,\quad
2\qN\a_{[Nt]-1}^\sN = \cO(N^{-3/2}),
$$
and invoking~\refp{asymptotic}, we get, for $k\in\cE$,
\begin{align}
\P\{T_{[Nt]}^\sN=k,S_{[Nt]}<0\}
&\sim
2\sqrt{\frac{2}{\pi}}\,\frac{\rho^-}{N^2} \sum_{i\in\cE:\atop k\le i\le [Nt]-2} \f\Big(\frac iN\Big)
\nonumber\\
&\hphantom{=\,}
+\frac{2}{\pi N^3}\sum_{i\in\cE:\atop k\le i\le [Nt]-2} \f\Big(\frac iN\Big)
\sum_{j\in\cE:\atop j\ge [Nt]-i-1} \f\Big(\frac jN\Big).
\label{asymptotic-intermediate1}
\end{align}
We can easily see that
\beq
\sum_{i\in\cE:\atop k\le i\le [Nt]-2} \f\Big(\frac iN\Big)
\sim \frac N2 \int_{k/N}^t \f(z)\,dz
\label{asymptotic-intermediate2}
\eeq
and
$$
\sum_{j\in\cE:\atop j\ge [Nt]-i-1} \f\Big(\frac jN\Big)
\sim \frac N2 \int_{([Nt]-i-1)/N}^{\infty} \f(v)\,dv
\sim \frac N2 \int_{t-i/N}^{\infty} \f(v)\,dv.
$$
Moreover,
\begin{align}
\sum_{i\in\cE:\atop k\le i\le [Nt]-2} \f\Big(\frac iN\Big)
\sum_{j\in\cE:\atop j\ge [Nt]-i-1} \f\Big(\frac jN\Big)
&\sim
\frac N2
\sum_{i\in\cE:\atop k\le i\le [Nt]-2} \f\Big(\frac iN\Big)
\int_{t-i/N}^{\infty} \f(v)\,dv
\nonumber\\
&\sim
\frac{N^2}{4} \int_{k/N}^t \f(u)\,du\int_{t-u}^{\infty} \f(v)\,dv.
\label{asymptotic-intermediate3}
\end{align}
Put now
$$
\chi^\pm(s,t)=\frac{2\rho^\pm}{\sqrt{2\pi}} \int_s^t \f(z)\,dz
+\frac{1}{2\pi} \int_s^t \f(u)\,du\int_{t-u}^{\infty} \f(v)\,dv.
$$
In the light of~\refp{asymptotic-intermediate1},
\refp{asymptotic-intermediate2} and~\refp{asymptotic-intermediate3},
we derive
$$
\P\{T_{[Nt]}^\sN=k,S_{[Nt]}<0\} \underset{N\to \infty}{\sim} \frac 1N\,\chi^-(k/N,t)
$$
and finally
$$
\P\{s_1< \frac 1N\,T_{[Nt]}^\sN\le s_2,B_t^\sN<0\}
\underset{N\to \infty}{\sim} \frac 1N \sum_{k\in\cE:\atop [Ns_1]<k\le [Ns_2]} \chi^-(k/N,t)
\underset{N\to \infty}{\rightarrow} \frac12 \int_{s_1}^{s_2} \chi^-(s,t)\,ds.
$$
Since $\frac 1N\,T_{[Nt]}^\sN \ind_{(-\infty,0)}(B_t^\sN)
\underset{N\to \infty}{\rightarrow} \bT_t \ind_{(-\infty,0)}(B_t)$, we deduce that
$$
\P\{s_1\le \bT_t \le s_2,B_t<0\}=\frac 12 \int_{s_1}^{s_2}\chi^-(s,t)\,ds
\quad\mbox{and}\quad \P\{\bT_t\in ds,B_t<0\}/ds=\frac 12\,\chi^-(s,t).
$$
We can prove in a quite similar way that
$$
\P\{\bT_t\in ds,B_t>0\}/ds=\frac 12\,\chi^+(t-s,t).
$$
This last result can be deduced also from the previous one by invoking
a duality argument related to drifted Brownian motion. Indeed, setting
more precisely $B_t=B_t^\r$ and $\bT_t=\bT_t^\r$, it can be easily seen that
$(B_t^\r,\bT_t^\r)$ and $(-B_t^{-\r},t-\bT_t^{-\r})$ have the same distributions.
This explains why the term $\chi^-(s,t)$ is changed into $\chi^+(t-s,t)$.
The proof of Theorem~\ref{limit-conditioned} is finished.
\fin

For $\r=0$, we retrieve the well-known results for
standard Brownian motion: for $s\in (0,t)$,
$$
\P\{\bT_t\in ds,B_t<0\}/ds=\frac{1}{\pi t}\sqrt{\frac{t-s}{s}}\quad\mbox{and}\quad
\P\{\bT_t\in ds,B_t>0\}/ds=\frac{1}{\pi t}\sqrt{\frac{s}{t-s}}.
$$

%%%%%%%%%%%%%%%%%%%%%%%%%%%%%%%%%%%%%%%%%%%%%%%%%%%%%%%%%%%%%%%%%%%%%%%%%%%%
%%%%%%%%%%%%%%%%%%%%%%%%%%%%%%%%%%%%%%%%%%%%%%%%%%%%%%%%%%%%%%%%%%%%%%%%%%%%
%%%%%%%%%%%%%%%%%%%%%%%%%%%%%%%%%%%%%%%%%%%%%%%%%%%%%%%%%%%%%%%%%%%%%%%%%%%%
\section{Appendix}

\demlem
The main idea is roughly speaking that, referring to Stirling formula for $a_i$,
$$
a_i=\frac{1}{i+2}\binom{i}{i/2} \underset{i\to \infty \atop i\in\cE}{\sim}
\frac{2^{i+1/2}}{\sqrt{\pi}\, i^{3/2}}\quad\mbox{and}\quad
(\pN\qN)^{i/2}=\bigg[\frac14\Big(1-\frac{\r^2}{N}\Big)\bigg]^{i/2}
\underset{i,N\to \infty\atop i/N^2\to 0}{\sim} \frac{1}{2^i}\,e^{-\r^2i/(2N)}.
$$
Then, recalling that $\a_i^\sN= a_i(\pN\qN)^{i/2}$,
$$
\a_i^\sN \underset{i,N\to \infty\atop i/N^2\to 0}{\sim}
\sqrt{\frac{2}{\pi}} \frac{e^{-\r^2i/(2N)}}{i^{3/2}}
$$
and next
$$
\sum_{i\in\cE:\atop i\ge k} \a_i^\sN
\underset{k,N\to \infty\atop c_1N\le k\le c_2N}{\sim} \sqrt{\frac{2}{\pi}}
\sum_{i\in\cE:\atop i\ge k} \frac{e^{-\r^2i/(2N)}}{i^{3/2}}
=\frac{\sqrt2}{\sqrt{\pi}\,N^{3/2}}
\sum_{i\in\cE:\atop i\ge k} \f(i/N)
\underset{k,N\to \infty\atop c_1N\le k\le c_2N}{\sim}
\sqrt{\frac{2}{\pi N}}\int_{k/N}^{\infty} \f(z)\,dz.
$$

For making the things more precise, let us write that
\begin{align*}
\Big(1-\frac{\r^2}{N}\Big)^{i/2}=\exp\bigg[\frac i2\,
\ln\!\Big(1-\frac{\r^2}{N}\Big)\bigg]
=e^{-\r^2i/(2N)}\,\exp\bigg[\frac i2\bigg(\ln\!\Big(1-
\frac{\r^2}{N}\Big)+\frac{\r^2}{N}\bigg)\bigg]
\end{align*}
from which we see that
$$
\Big(1-\frac{\r^2}{N}\Big)^{i/2}
\underset{i,N\to \infty\atop i/N^2\to 0}{\sim}e^{-\r^2i/(2N)}.
$$
Pick $\e>0$. There exists an integer $N_0$ such that, for any $N\ge N_0$,
for any $i,k\in\cE$ such that $c_1N\le k\le i \le N^{3/2}$,
$$
\a_i^\sN\ge (1-\e)\,\frac{\sqrt2}{\sqrt{\pi}\,{i^{3/2}}}\,e^{-\r^2i/(2N)}
$$
and for any $i,k\in\cE$ such that $c_1N\le k\le i\wedge (c_2N)$,
$$
\a_i^\sN\le (1+\e)\,\frac{\sqrt2}{\sqrt{\pi}\,{i^{3/2}}}\,e^{-\r^2i/(2N)}.
$$
Then, for any $N\ge N_0$ and any $k$ such that $c_1N\le k\le c_2N$,
\beq
(1-\e)\,\frac{\sqrt2}{\sqrt{\pi}\,{N^{3/2}}}
\sum_{i\in\cE:\atop k\le i\le N^{3/2}} \f\Big(\frac iN\Big)
\le \sum_{i\in\cE:\atop i\ge k} \a_i^\sN
\le (1+\e)\,\frac{\sqrt2}{\sqrt{\pi}\,{N^{3/2}}}
\sum_{i\in\cE:\atop i\ge k} \f\Big(\frac iN\Big).
\label{bounds1}
\eeq
Remarking that $\f$ is decreasing, we have
$$
\frac N2 \int_{i/N}^{(i+2)/N} \f(z)\,dz \le\f\Big(\frac iN\Big)
\le \frac N2 \int_{(i-2)/N}^{i/N} \f(z)\,dz.
$$
Then
$$
\frac N2 \int_{k/N}^{\infty} \f(z)\,dz \le
\sum_{i\in\cE:\atop i\ge k}\f\Big(\frac iN\Big)
\le \f\Big(\frac kN\Big)+\frac N2 \int_{k/N}^{\infty} \f(z)\,dz
$$
which shows, since we plainly have for $c_1N\le k\le c_2N$,
$\f(k/N)=o(N\int_{k/N}^{\infty} \f(z)\,dz)$, that
\beq
\sum_{i\in\cE:\atop i\ge k} \f\Big(\frac iN\Big)
\underset{k,N\to \infty\atop c_1N\le k\le c_2N}{\sim}
\frac N2 \int_{k/N}^{\infty} \f(z)\,dz.
\label{bounds2}
\eeq
On the other hand,
\beq
\sum_{i\in\cE:\atop k\le i\le N^{3/2}} \f\Big(\frac iN\Big)
=\sum_{i\in\cE:\atop i\ge k} \f\Big(\frac iN\Big)
-\sum_{i\in\cE:\atop i> N^{3/2}} \f\Big(\frac iN\Big).
\label{bounds3}
\eeq
The last sum in~\refp{bounds3} can be estimated as follows:
\begin{align}\label{bounds4}
\sum_{i\in\cE:\atop i> N^{3/2}}\f\Big(\frac iN\Big)
&\le
\frac N2 \int_{\sqrt N}^{\infty} \f(z)\,dz +\f\Big(\frac{[N^{3/2}]}{N}\Big)
+\f\Big(\frac{[N^{3/2}]+1}{N}\Big)
\nonumber\\
&\underset{N\to \infty}{\sim}
\frac N2 \int_{\sqrt N}^{\infty} \f(z)\,dz
\underset{N\to \infty}{=}o(N).
\end{align}
Putting~\refp{bounds2} and~\refp{bounds4} into~\refp{bounds3}, we obtain that
\beq
\sum_{i\in\cE:\atop k\le i\le N^{3/2}} \f\Big(\frac iN\Big)
\underset{N\to \infty}{\sim}
\frac N2 \int_{k/N}^{\infty} \f(z)\,dz
\label{bounds5}
\eeq
and next putting~\refp{bounds2} and~\refp{bounds5} into~\refp{bounds1},
we find, for $N$ large enough and \mbox{$c_1N\!\le k\le\! c_2N$,} that
$$
(1-2\e)\,\frac{\sqrt2}{\sqrt{\pi}\,N} \int_{k/N}^{\infty} \f(z)\,dz
\le \sum_{i\in\cE:\atop i\ge k} \a_i^\sN
\le (1+2\e)\,\frac{\sqrt2}{\sqrt{\pi}\,N} \int_{k/N}^{\infty} \f(z)\,dz
$$
from which we finally deduce~\refp{asymptotic}.
\fin

%%%%%%%%%%%%%%%%%%%%%%%%%%%%%%%%%%%%%%%%%%%%%%%%%%%%%%%%%%%%%%%%%%%%%%%%%%%%
%%%%%%%%%%%%%%%%%%%%%%%%% REFERENCES %%%%%%%%%%%%%%%%%%%%%%%%%%%%%%%%%%%%%%%
%%%%%%%%%%%%%%%%%%%%%%%%%%%%%%%%%%%%%%%%%%%%%%%%%%%%%%%%%%%%%%%%%%%%%%%%%%%%

\end{document}